\documentclass[12pt]{amsart}

\usepackage{amsmath,amssymb,amsthm}
\usepackage{graphicx}
\usepackage[hypertexnames=false]{hyperref}
\hypersetup{pdftitle={The error in a smooth weighted prime number formula and zero-free regions for the Riemann zeta function}}
\usepackage{cleveref}
\usepackage{fancyhdr}
\usepackage{xcolor}

\usepackage[margin=1.5in]{geometry}

\theoremstyle{plain}
\newtheorem{theorem}{Theorem}[section]
\newtheorem{lemma}[theorem]{Lemma}
\newtheorem{proposition}[theorem]{Proposition}
\newtheorem{corollary}[theorem]{Corollary}

\theoremstyle{definition}
\newtheorem{definition}[theorem]{Definition}

\theoremstyle{remark}

\theoremstyle{plain}
\newtheorem{thmA}{Theorem}

\newcommand{\eps}{\varepsilon}

\title{%
    The error in a smooth weighted prime number formula\\
    and zero-free regions for the Riemann zeta function
}
\date{}

\author{Songlin Han}

\begin{document}
\sloppy
\hfuzz=0.1pt
\vfuzz=1pt

\maketitle

\begin{abstract}
    We study the error bound for a smooth weighted prime number formula, $\sum_{n}(\Lambda(n)-1)\exp(-n/x)$, and its implication to the zero-free region for the Riemann zeta function by improving the method of Pintz.
    We also give an application to the average number of smooth weighted Goldbach representations and generalize to smooth weighted average $k$-Goldbach representations.
\end{abstract}

\section{Introduction}
Let $\Lambda(n)$ denote the von Mangoldt function, defined by 
\begin{equation*}
    \Lambda(n)=
        \begin{cases}
            \log p & \text{if } n=p^m \text{ for some prime $p$ and }m\ge 1,\\
            0 & \text{otherwise.}
        \end{cases}
\end{equation*}
This function plays a central role in analytic number theory, particularly in studying the distribution of prime numbers. The classical prime number theorem (PNT) asserts the asymptotic formula
\begin{equation}
    \sum_{n \leq x} \Lambda(n) \sim x.
\end{equation}
In other words, $\Lambda(n)$ averages out to $1$ on the scale of natural numbers. This naturally raises the question: how does $\Lambda(n)$ deviate from this average?

To explore such deviations, one may consider the difference $\Lambda(n) - 1$. This idea has appeared in various forms, including in a formula proposed by Hardy-Littlewood \cite{HardyLittlewood1916} in the early 20th century. They introduced the function
\begin{equation*}
    \Delta(z) := \sum_{n=1}^\infty \left(\Lambda(n) - 1\right) e^{-\frac{n}{z}},
\end{equation*}
where $z$ is a complex number with $\Re(z) > 0$. This allows one to focus on finer aspects of the oscillation of $\Lambda(n)$.

To analyze $\Delta(x)$ for real $x > 0$, it is convenient to define the following two functions:
\begin{equation*}
    \Psi(x) := \sum_{n=1}^\infty \Lambda(n) e^{-\frac{n}{x}}
    \quad\text{and}\quad
    I(x) := \sum_{n=1}^\infty e^{-\frac{n}{x}}.
\end{equation*}
Then we have the identity
\begin{equation}
    \Delta(x) = \Psi(x) - I(x).
\end{equation}
Here, $\Psi(x)$ captures a weighted contribution of primes, while $I(x) = \frac{1}{e^{1/x}-1} = x - 1/2 + O(1/x)$ acts as a smooth baseline---their difference isolates the irregularity in the distribution of prime powers.

An explicit formula for $\Psi(x)$, originally derived by Hardy and Littlewood \cite{HardyLittlewood1916} (see also \cite[12.1.1 Exercise 8 (c)]{MontgomeryVaughan2007}), reveals deeper connections with the non-trivial zeros of the Riemann zeta function $\zeta(s)$. Specifically,
\begin{equation}\label{explicit_formula}
    \Psi(x) = x - \sum_{\rho}\Gamma(\rho)x^{\rho} -\log 2\pi + O\left(\frac{1}{x}\right),
\end{equation}
where $\rho$ runs over all non-trivial zeros of $\zeta(s)$, and $\Gamma$ is the gamma function. From this, they obtained the following conditional result under the Riemann Hypothesis (RH).

\begin{thmA}[Theorem 2.24, \cite{HardyLittlewood1916}]
    \label{thm:HL}
    If the Riemann hypothesis is true, then $\Delta(x) = O\left(\sqrt{x}\right)$ as $x$ tends to $\infty$.
\end{thmA}

Instead of assuming RH, we relate arbitrary zero-free regions (ZFR) for the Riemann zeta function to the exponential decreasing of the error term $\Delta(x)$. The unsmoothed analogue of this connection is fully studied: the equivalence between ZFRs of $\zeta$ and the size of
\begin{equation*}
    R(x) := \sum_{n\le x}\Lambda(n) - x
\end{equation*}
(the remainder term in the classical PNT formula) was established in a sequence of works by Ingham, Tur\'an, Sta\'s, and Pintz. We conclude their result in the form given in \cite[Section 1]{BillingtonChengSchettlerSuriajaya2025}.

\begin{thmA}[Classical PNT--ZFR; Ingham \cite{Ingham1990}, Tur\'an \cite{Turan1950}, Sta\'s \cite{Stas1961}, Pintz \cite{Pintz1980I,Pintz1980II,Pintz1984}]
    \label{thm:PNTZFR}
    Fix $0 < \varepsilon < 1$, and let $\eta:[0,\infty)\to(0,1/2]$ be continuous and decreasing.
    \begin{itemize}
        \item If $\zeta(\sigma+it)\neq 0$ in the region $\sigma > 1 - \eta(|t|)$, then
        \begin{equation*}
            R(x) \ll x\exp\bigl(-(1-\varepsilon)\omega(x)\bigr), \qquad \omega(x) := \inf_{u\ge 1}\bigl(\eta(u)\log x + \log u\bigr).
        \end{equation*}
        \item If $R(x) \ll x\exp\bigl(-(1+\varepsilon)\varpi(x)\bigr)$ with
        \begin{equation*}
            \varpi(x) := \inf_{u\ge 0}\bigl(\eta(u)\log x + u\bigr),
        \end{equation*}
        and $\lim_{u\to\infty}\eta'(u) = 0$, then for all sufficiently large $|t|$, $\zeta(\sigma+it)\neq 0$ in the region $\sigma > 1 - \eta(\log|t|)$.
    \end{itemize}
\end{thmA}

Our goal is to prove the analogue of Theorem~\ref{thm:PNTZFR} with the classical remainder $R(x)$ replaced by the exponentially smoothed remainder $\Delta(x)$. The smoothed setting leads to a single unified localization quantity $\omega_\eta$ that uses $\log|\Gamma(it)|$ in place of the separate $\log u$ and $u$ of Theorem~\ref{thm:PNTZFR}, reflecting the presence of $\Gamma(\rho)$ in the smoothed explicit formula \eqref{explicit_formula}. The converse direction of our theorem requires additional conditions on $\eta$, which we formalize as follows.

\begin{definition}[Admissibility]\label{def:admissible}
    A function $\eta:[2,\infty)\to(0,1/2]$ is called admissible if, for all sufficiently large~$t$:
    \begin{enumerate}
        \item[\textup{(H1)}] $\eta$ is $C^2$ and strictly decreasing: $\eta'(t) < 0$.
        \item[\textup{(H2)}] $\eta$ is strictly convex: $\eta''(t) > 0$.
        \item[\textup{(H3)}] Slow variation: $\displaystyle\lim_{t\to\infty}\frac{t\eta'(t)}{\eta(t)} = 0$.
        \item[\textup{(H4)}] $\displaystyle\lim_{t\to\infty}\eta(t) = 0$.
    \end{enumerate}
\end{definition}

\begin{theorem}[Smooth weighted PNT and ZFR]
    \label{thm:main1}
    Let $0 < \varepsilon < 1$ be fixed. Let $\eta: [2, \infty) \to (0, 1/2]$ be a continuously differentiable, monotonically decreasing function. Define
    \begin{equation} \label{eq:omega_eta}
        \omega_{\eta}(x) := \inf_{t \ge 2} \left( \eta(t)\log x - \log |\Gamma(it)| \right).
    \end{equation}
    The error term $\Delta(x)$ and the zero-free region for the Riemann zeta function are related as follows.
    \begin{itemize}
        \item  If $\zeta(\sigma+it)\neq 0$ in the region $\sigma > 1 - \eta(|t|)$ for $|t| \ge 2$, then 
        \begin{equation*}
            \Delta(x) \ll x\exp\left(-(1-\varepsilon)\omega_{\eta}(x)\right)
        \end{equation*}
        as $x$ tends to $\infty$.

        \item  Suppose $\eta$ is admissible in the sense of Definition~\ref{def:admissible}, and
        \begin{equation*}
            \Delta(x) \ll x\exp\left(-(1+\varepsilon)\omega_{\eta}(x)\right).
        \end{equation*}
        Then for every $\delta > 0$ there exists $t_0 = t_0(\delta) \ge 2$ such that $\zeta(\sigma+it)\neq 0$ in the region 
        \begin{equation*}
            \sigma > 1 - (1-\delta)\eta(|t|)
        \end{equation*}
        for all $|t| \ge t_0$.
    \end{itemize}
\end{theorem}

For a fixed half-plane zero-free region, as in the quasi-RH, we have the following equivalence: for any fixed $0<\theta<1/2$, all non-trivial zeros satisfy $\beta\le1-\theta$ if and only if $\Delta(x)=O(x^{1-\theta})$ as $x\to\infty$. The forward implication follows from \eqref{explicit_formula}, and the converse from the Mellin transform.
We therefore focus on zero-free regions whose boundary approaches the line $\sigma=1$, with $\eta$ satisfying the admissibility conditions in the converse direction.

\bigskip

Let us consider a different type of convolution of $\Lambda$. For an integer $k \ge 2$, define a $k$-times additive convolution of the von Mangoldt function
\begin{equation}
    \psi_k(n) := \sum_{m_1 + m_2 + \dots + m_k = n} \Lambda(m_1)\Lambda(m_2)\cdots\Lambda(m_k).
\end{equation}
When $k=2$, this function weights the number of ways to express an integer $n$ as a sum of two prime powers. The cumulative sum $G(x) := \sum_{n\le x}\psi_2(n)$ is often referred to as the average number of Goldbach representations. This average number of Goldbach representations was studied by Fujii \cite{Fujii1991}, Granville \cite{Granville2007}, and also \cite{FriedlanderGoldston1997,LanguascoZaccagnini2012a,LanguascoZaccagnini2012b,Nguyen2024}.

Recently, the connection between ZFRs and the average number of Goldbach representations has been explored using the circle method \cite{BillingtonChengSchettlerSuriajaya2025}, relying on contour integral identities over the Fej\'er kernel \cite{GoldstonSuriajaya2023}. In \cite{BillingtonChengSchettlerSuriajaya2025}, an asymptotic formula for the following smoothed version of the number of Goldbach representations is given.

\begin{thmA}[Theorem 3, \cite{BillingtonChengSchettlerSuriajaya2025}]
    \label{thm:GR}
    Suppose that $\zeta(\sigma+it)\neq 0$ in $\sigma>1-\eta(|t|)$ for some continuous decreasing function $\eta(u)$ such that $0<\eta(u)\le 1/2$. Then for $N\ge 3$,
    \begin{equation}
        \sum_{n=1}^{\infty}\psi_2(n)e^{-n/N} = N^2 + O\left(N^{2-\eta(\log N)}\right).
    \end{equation}
\end{thmA}

The authors in \cite{BillingtonChengSchettlerSuriajaya2025} suggested that one can obtain a converse to this theorem by the method of Pintz \cite{Pintz1980II,Pintz1984}. In this paper, we establish this converse. This observation allows us to bridge the smooth weighted Goldbach representations to Theorem~\ref{thm:main1}. As a result, the precise equivalence naturally generalizes to smooth $k$-Goldbach representations.

\begin{theorem}[Smooth weighted $k$-Goldbach representations and ZFR]
    \label{thm:main2}
    Let $k \ge 2$ be an integer, and let $0 < \varepsilon < 1$. Let $\eta(t)$ and $\omega_\eta(x)$ be defined exactly as in Theorem~\ref{thm:main1}. Define the smooth weighted average $k$-Goldbach sum:
    \begin{equation*}
        G_k(x) := \sum_{n=1}^{\infty} \psi_k(n)e^{-n/x}.
    \end{equation*}
    \begin{itemize}
        \item If $\zeta(\sigma+it)\neq 0$ in the region $\sigma > 1-\eta(|t|)$ for $|t| \ge 2$, then 
        \begin{equation*}
            G_k(x) = x^k + O_k\left(x^k \exp\left(-(1-\varepsilon)\omega_{\eta}(x)\right)\right).
        \end{equation*}
        \item Suppose $\eta$ is admissible (Definition~\ref{def:admissible}), and the above asymptotic holds with exponent $(1+\varepsilon)$, i.e.,
        \begin{equation*}
            G_k(x) = x^k + O_k\left(x^k \exp\left(-(1+\varepsilon)\omega_{\eta}(x)\right)\right).
        \end{equation*}
        Then for every $\delta > 0$ there exists $t_0 = t_0(\delta) \ge 2$ such that $\zeta(\sigma+it)\neq 0$ in the region $\sigma > 1-(1-\delta)\eta(|t|)$ for all $|t| \ge t_0$.
    \end{itemize}
\end{theorem}

\bigskip
The paper is organized as follows. Section~\ref{sec:prelim} shows the preliminary results which will be used throughout. In Section~\ref{sec:proofs} we prove Theorems~\ref{thm:main1} and~\ref{thm:main2}. Section~\ref{sec:upper_bound} establishes the auxiliary bound on $W(x) := \sum_{0<\gamma\le x}|\Gamma(\rho+1)|x^\beta/\gamma$ used in the upper-bound direction of Theorem~\ref{thm:main1}.

\section{Preliminaries}
\label{sec:prelim}

\subsection{Notation and Stirling asymptotics}
\label{subsec:notation}

Throughout this paper, $\rho = \beta + i\gamma$ with $\gamma,\beta\in\mathbb{R}$ denotes a non-trivial zero of $\zeta(s)$. We use standard Landau ($O, o$) and Vinogradov ($\ll, \gg, \sim$) notations for asymptotic estimates. Unless otherwise specified, the summation symbol $\sum_n$ is taken over all positive integers, i.e.,
\begin{equation*}
    \sum_{n} = \sum_{n=1}^{\infty}.
\end{equation*}
Similarly, $\sum_{\rho}$ denotes a sum over all non-trivial zeros of $\zeta$ in this paper, implying the truncation $\lim_{T\to\infty}\sum_{|\gamma|\le T}$ when necessary. We write $N(T) := \#\{\rho : 0 < \gamma \le T\}$ for the zero-counting function, and recall the Riemann--von Mangoldt bounds $N(T) \ll T\log T$ and $N(T+1) - N(T) \ll \log T$.

We describe some facts of the gamma function from \cite[Chapter~2]{Remmert1998}.
Stirling's formula plays a very important role in this paper. We use it in the form
\begin{equation}\label{eq:stirling_master}
    |\Gamma(\sigma+it)| = \sqrt{2\pi}|t|^{\sigma-1/2}e^{-\pi|t|/2}\left(1 + O\left(\frac{1}{|t|}\right)\right)
\end{equation}
for $0 \le \sigma \le 1$ as $|t|$ tends to $\infty$, and we write $g(t) := -\log|\Gamma(it)|$ for convenience.
For real $t > 0$, Schwarz reflection gives $\overline{\Gamma(it)} = \Gamma(-it)$, so $|\Gamma(it)|^2 = \Gamma(it)\Gamma(-it)$. Using the functional equation $\Gamma(1-it) = (-it)\Gamma(-it)$, followed by $\Gamma(s)\Gamma(1-s) = \pi/\sin(\pi s)$ at $s = it$, and $\sin(\pi it) = i\sinh(\pi t)$, we obtain
\begin{equation}\label{eq:gamma_it_identity}
    |\Gamma(it)|^2 = \frac{\Gamma(it)\Gamma(1-it)}{-it} = \frac{\pi}{-it\cdot i\sinh(\pi t)} = \frac{\pi}{t\sinh(\pi t)}.
\end{equation}
Taking $-\tfrac{1}{2}\log$ of both sides gives the exact formula
\begin{equation}\label{eq:g_exact}
    g(t) = \tfrac{1}{2}\log t + \tfrac{1}{2}\log\sinh(\pi t) - \tfrac{1}{2}\log\pi, \qquad g'(t) = \frac{1}{2t} + \frac{\pi}{2}\coth(\pi t),
\end{equation}
valid for all $t > 0$. In particular, since $\coth(\pi t) > 1$ for $t > 0$, we have the useful uniform bound
\begin{equation}\label{eq:g_prime_lower}
    g'(t) > \pi/2 \qquad\text{for all }t > 0.
\end{equation}
To get asymptotic estimates, we write $\sinh(\pi t) = \tfrac{1}{2}e^{\pi t}(1 - e^{-2\pi t})$, thus
\begin{equation*}
    \tfrac{1}{2}\log\sinh(\pi t) = \tfrac{\pi t}{2} - \tfrac{1}{2}\log 2 + \tfrac{1}{2}\log(1 - e^{-2\pi t}) = \tfrac{\pi t}{2} - \tfrac{1}{2}\log 2 + O(e^{-2\pi t}).
\end{equation*}
Combining with \eqref{eq:g_exact} and $\tfrac{1}{2}\log 2 + \tfrac{1}{2}\log\pi = \tfrac{1}{2}\log(2\pi)$, we have
\begin{equation}\label{eq:g_asymptotics}
    g(t) = \tfrac{\pi}{2}t + \tfrac{1}{2}\log t - \tfrac{1}{2}\log(2\pi) + O(e^{-2\pi t}), \qquad g'(t) = \tfrac{\pi}{2} + \tfrac{1}{2t} + O(e^{-2\pi t}),
\end{equation}
as $t $ tends to $\infty$.

\subsection{Localization of the infimum}

For an admissible $\eta$ and $\gamma > 0$, define
\begin{equation}\label{eq:L_gamma}
    L=L_\gamma := \frac{\pi}{2|\eta'(\gamma)|}, \qquad f(t)=f_\gamma(t) := \eta(t) L_\gamma + g(t),
\end{equation}
\begin{equation}\label{eq:AB}
    A_\gamma := \eta(\gamma) L_\gamma = \frac{\pi\eta(\gamma)}{2|\eta'(\gamma)|}, \qquad B_\gamma := g(\gamma) = -\log|\Gamma(i\gamma)|.
\end{equation}
By (H3),
\begin{equation}\label{eq:ratio}
    \frac{A_\gamma}{\gamma} = \frac{\pi}{2}\cdot\frac{\eta(\gamma)}{\gamma|\eta'(\gamma)|} \longrightarrow +\infty,
\end{equation}
so $A_\gamma \to \infty$ as $\gamma \to \infty$. Since $B_\gamma \sim \frac{\pi}{2}\gamma$ by \eqref{eq:g_asymptotics}, we have $B_\gamma = o(A_\gamma)$.

\begin{lemma}[Localization]\label{lem:localization}
    For $\eta$ admissible and $\gamma$ sufficiently large,
    \begin{equation}\label{eq:localization}
        \omega_\eta(e^{L_\gamma}) = A_\gamma + B_\gamma - o(A_\gamma).
    \end{equation}
\end{lemma}

\begin{proof}
    Setting $x = e^L$ in \eqref{eq:omega_eta} and using $g(t) = -\log|\Gamma(it)|$,
    \begin{equation*}
        \omega_\eta(e^L) = \inf_{t\ge 2}\bigl(\eta(t)L + g(t)\bigr) = \inf_{t\ge 2} f(t).
    \end{equation*}
    Evaluating at $t = \gamma$, $f(\gamma) = \eta(\gamma)L + g(\gamma) = A_\gamma + B_\gamma$ by \eqref{eq:AB}, which gives the upper bound $\omega_\eta(e^L) \le A_\gamma + B_\gamma$. We prove the matching lower bound
    \begin{equation*}
        \inf_{t \ge 2} f(t) \ge A_\gamma + B_\gamma - o(A_\gamma).
    \end{equation*}

    We first consider the case $t\in[2,\gamma]$.
    Since $\eta$ is decreasing, $\eta(t) \ge \eta(\gamma)$ for $t \in [2,\gamma]$, so $\eta(t)L \ge A_\gamma$. By \eqref{eq:g_prime_lower}, $g'(t) > 0$ for all $t > 0$, so $g$ is strictly increasing on $(0,\infty)$, in particular $g(t) \ge g(2)$ on $[2,\gamma]$. Combining these, we have
    \begin{equation*}
        f(t) \ge A_\gamma + g(2) = A_\gamma + B_\gamma - (B_\gamma - g(2)), \qquad t \in [2,\gamma].
    \end{equation*}
    By \eqref{eq:g_asymptotics} we have $B_\gamma \sim \tfrac{\pi}{2}\gamma$, and by \eqref{eq:ratio}, $\gamma = o(A_\gamma)$. Hence $B_\gamma - g(2) = o(A_\gamma)$, giving $f(t) \ge A_\gamma + B_\gamma - o(A_\gamma)$ on $[2,\gamma]$.

    \medskip\noindent
    We now turn to the case $t \in [\gamma,\infty)$.
    We show $f$ is strictly increasing on $[\gamma,\infty)$, it then follows that $f(t) \ge f(\gamma) = A_\gamma + B_\gamma$ for all $t \ge \gamma$.

    By (H2), $\eta'$ is strictly increasing, so for $t \ge \gamma$,
    \begin{equation*}
        \eta'(t)L \ge \eta'(\gamma)L = -\tfrac{\pi}{2},
    \end{equation*}
    using the choice $L = \pi/(2|\eta'(\gamma)|)$. By \eqref{eq:g_prime_lower}, $g'(t) > \pi/2$ for all $t > 0$. Therefore
    \begin{equation*}
        f'(t) = \eta'(t)L + g'(t) > -\tfrac{\pi}{2} + \tfrac{\pi}{2} = 0, \qquad t \ge \gamma,
    \end{equation*}
    so $f$ is strictly increasing on $[\gamma,\infty)$, and $f(t) \ge A_\gamma + B_\gamma$ in this case.

    \medskip
    Combining the two cases we have $\inf_{t\ge 2} f(t) \ge A_\gamma + B_\gamma - o(A_\gamma)$.
\end{proof}

\begin{lemma}[Consequence of slow variation]\label{lem:slow_var}
    Let $\eta$ be admissible. Then, as $t\to\infty$,
    \begin{equation*}
        |\log\eta(t)|=o(\log t), \qquad \eta(t)=t^{-o(1)}.
    \end{equation*}
    In particular, $t^q\eta(t)\to\infty$ for every fixed $q>0$.
\end{lemma}

\begin{proof}
    Fix $\kappa>0$. By (H3), there exists $T_\kappa>1$ such that
    \begin{equation*}
        \left|\frac{\eta'(u)}{\eta(u)}\right|\le\frac{\kappa}{u}
        \qquad (u\ge T_\kappa).
    \end{equation*}
    Since $\eta$ is positive, integration gives, for $t\ge T_\kappa$,
    \begin{equation*}
        |\log\eta(t)-\log\eta(T_\kappa)|
        \le \kappa\log(t/T_\kappa).
    \end{equation*}
    Hence
    \begin{equation*}
        \limsup_{t\to\infty}\frac{|\log\eta(t)|}{\log t}\le\kappa.
    \end{equation*}
    As $\kappa>0$ is arbitrary, this proves the first assertion.
    Since $0<\eta(t)\le1/2$, it also gives $\eta(t)=t^{-o(1)}$.
    Finally, for every fixed $q>0$,
    $\log(t^q\eta(t))=q\log t+o(\log t)\to\infty$.
\end{proof}

\subsection{Tur\'an's power-sum method}
\label{subsec:turan}

The converse of Theorem~\ref{thm:main1} relies on the following classical power-sum inequality, which we quote in the form given in Montgomery's lectures \cite{Montgomery1994}.

\begin{thmA}[Tur\'an's Second Main Theorem; {\cite[p.~94]{Montgomery1994}}]\label{thm:turan}
    Let $s_\nu=\sum_{n=1}^N b_n z_n^\nu$, where $b_n,z_n\in\mathbb{C}$ and $\nu$ is a non-negative integer. Suppose $1=|z_1|\ge|z_2|\ge\cdots\ge|z_N|$. Then for every non-negative integer $M$, there exists $\nu$ with $M+1\le\nu\le M+N$ such that
    \begin{equation}\label{eq:turan_discrete}
        |s_\nu|\ge 2\left(\frac{N}{8e(M+N)}\right)^{N}\min_{1\le j\le N}\left|\sum_{n=1}^j b_n\right|.
    \end{equation}
\end{thmA}

We will use the following corollary in the next section.

\begin{corollary}\label{cor:turan_cont}
    Let $W(\mu)=\sum_{i=1}^N b_i e^{\lambda_i\mu}$ with $\Re(\lambda_1)\ge\cdots\ge\Re(\lambda_N)$ and $\Re(\lambda_1)\ge 0$. Then for any $a\ge 0$ and $d>0$,
    \begin{equation}\label{eq:turan_cont}
        \max_{a\le\mu\le a+d}|W(\mu)|\ge 2e^{a\Re(\lambda_1)}\left(\frac{d}{16e(a+d)}\right)^{N}\min_{1\le\ell\le N}\left|\sum_{i=1}^\ell b_i\right|.
    \end{equation}
\end{corollary}

\begin{proof}
    Set $h:=d/(N+1)$ and $z_i:=e^{(\lambda_i-\lambda_1)h}$, so that $|z_1|=1\ge|z_i|$. Let $M:=\lfloor a/h\rfloor$. For $\nu\in[M+1,M+N]$, set $\mu_\nu:=\nu h\in(a,a+d)$ and define
    \begin{equation*}
        \tilde s_\nu:=e^{-\lambda_1\mu_\nu}W(\mu_\nu)=\sum_i b_i z_i^\nu.
    \end{equation*}
    Theorem~\ref{thm:turan} yields $\nu\in[M+1,M+N]$ with
    \begin{equation*}
        |\tilde s_\nu|\ge 2\left(\frac{N}{8e(M+N)}\right)^N\min_\ell\left|\sum_{i=1}^\ell b_i\right|.
    \end{equation*}
    Since $|W(\mu_\nu)|=|\tilde s_\nu|e^{\Re(\lambda_1)\mu_\nu}\ge|\tilde s_\nu|e^{a\Re(\lambda_1)}$ and
    \begin{equation*}
        \frac{N}{M+N}\ge\frac{Nd}{(N+1)(a+d)}\ge\frac{d}{2(a+d)},
    \end{equation*}
    this completes the proof.
\end{proof}

\subsection{Analysis of \texorpdfstring{$\Gamma$ and $H$}{Gamma and H} near a fixed zero}
\label{subsec:coeff}

Throughout this subsection, $\rho_j = \beta_j + i\gamma_j$ denotes a fixed non-trivial zero of $\zeta$ with $\gamma_j > 0$ sufficiently large, and $H_j := H(\rho_j)$, where $H(s) := \Gamma'(s)/\Gamma(s)$ is the digamma function.
Using the Stirling formula we have the digamma asymptotic $H(s) = \log s + O(|s|^{-1})$, and
\begin{equation}\label{eq:digamma_asym}
    H_j = \sigma_j + i\nu_j, \qquad \sigma_j = \log\gamma_j + O(\gamma_j^{-1}), \qquad \nu_j = \frac{\pi}{2} + O(\gamma_j^{-1}).
\end{equation}
For any non-trivial zero $\rho = \beta + i\gamma$ of $\zeta$, we write $a_\rho := \beta - \beta_j$, $b_\rho := \gamma - \gamma_j$, and define
\begin{equation}\label{eq:tildec_def}
    \tilde c_\rho := \Gamma(\rho)e^{-H_j(\rho-\rho_j)}.
\end{equation}

    We first study the local behavior of $\tilde c_\rho$ near $\rho_j$.
    To do this, we need the following lemma.

\begin{lemma}[Local expansion of $\Gamma$]\label{lem:gamma_local}
    Let $w=u+iv$ with $u\ge0$ and $v\ge1$, and let $H=\Gamma'/\Gamma$ denote the
    digamma function.
    \begin{enumerate}
        \item[\textup{(a)}] For all $s$ with $|s-w|\le v/2$,
        \begin{equation}\label{eq:gamma_local_log}
            \left|\log\frac{\Gamma(s)e^{-H(w)(s-w)}}{\Gamma(w)}\right|
            \le \frac{3|s-w|^2}{v},
        \end{equation}
        where $\log$ denotes the branch that vanishes at $s=w$.
        \item[\textup{(b)}] If moreover $|s-w|^2\le v/12$, then
        \begin{equation}\label{eq:gamma_local_mult}
            \frac{\Gamma(s)e^{-H(w)(s-w)}}{\Gamma(w)} = 1+\theta,
            \qquad |\theta|\le\frac{4|s-w|^2}{v}.
        \end{equation}
        In particular, \eqref{eq:gamma_local_mult} controls both the modulus and
        the argument of the left-hand side:
        $\bigl||1+\theta|-1\bigr|\le|\theta|$ and
        $|\arg(1+\theta)|\le\tfrac{\pi}{2}|\theta|$.
    \end{enumerate}
\end{lemma}

\begin{proof}
    Since $u\ge0$, every pole $-k$ of $\Gamma$ (with $k\ge0$ an integer) satisfies $|w+k|=\sqrt{(u+k)^2+v^2}\ge v$, so the open disc $|s-w|<v$ contains no poles of $\Gamma$.
    As $\Gamma$ is zero-free and the disc is simply connected, there is a holomorphic branch of $\log\Gamma$ on this disc, whose Taylor expansion
    \begin{equation*}
        \log\Gamma(s)
        = \log\Gamma(w) + H(w)(s-w)
          + \sum_{m\ge2}\frac{H^{(m-1)}(w)}{m!}(s-w)^m
    \end{equation*}
    converges, in particular on $|s-w|\le v/2$.
    Writing
    \begin{equation}\label{eq:F_def}
        F(s) := \sum_{m\ge2}\frac{H^{(m-1)}(w)}{m!}(s-w)^m,
        \qquad\text{so that}\qquad
        e^{F(s)} = \frac{\Gamma(s)e^{-H(w)(s-w)}}{\Gamma(w)},
    \end{equation}
    the left-hand side of \eqref{eq:gamma_local_log} equals $|F(s)|$, and
    \begin{equation}\label{eq:F_tail}
        |F(s)| \le \sum_{m\ge2}\frac{|H^{(m-1)}(w)|}{m!}\,|s-w|^m.
    \end{equation}

    Differentiating the partial-fraction
    expansion
    $H(z) = -\gamma_{\mathrm{E}} - \frac{1}{z} + \sum_{k\ge1}\bigl(\frac{1}{k}-\frac{1}{z+k}\bigr)$ term by term, where $\gamma_{\mathrm{E}}$ is Euler's constant gives, for $m\ge2$,
    \begin{equation*}
        H^{(m-1)}(z) = (-1)^m(m-1)!\sum_{k=0}^\infty(z+k)^{-m}.
    \end{equation*}
    Since $u\ge0$, we have $|w+k|^2 = (u+k)^2+v^2 \ge k^2+v^2$ for every $k\ge0$.
    The function $t\mapsto(t^2+v^2)^{-m/2}$ is decreasing on $[0,\infty)$, so
    comparing the sum over $k\ge1$ with the integral and substituting $t=vx$,
    \begin{align*}
        \sum_{k=0}^\infty|w+k|^{-m}
        &\le v^{-m} + \int_0^\infty\frac{dt}{(t^2+v^2)^{m/2}} \\
        &= v^{-m} + v^{1-m}\int_0^\infty\frac{dx}{(1+x^2)^{m/2}} \\
        &\le v^{1-m}\Bigl(\frac1v+\frac{\pi}{2}\Bigr),
    \end{align*}
    where we used $(1+x^2)^{-m/2}\le(1+x^2)^{-1}$ for $m\ge2$. As $v\ge1$, the
    bracket is at most $1+\pi/2\le3$, hence
    \begin{equation}\label{eq:polygamma_bound}
        |H^{(m-1)}(w)| \le 3\,(m-1)!\,v^{1-m}
        \qquad (m\ge2).
    \end{equation}

    Inserting \eqref{eq:polygamma_bound} into \eqref{eq:F_tail} and writing
    $r:=|s-w|$ and $q:=r/v\le1/2$,
    \begin{equation*}
        |F(s)| \le 3v\sum_{m\ge2}\frac{q^m}{m}
        \le \frac{3v}{2}\sum_{m\ge2}q^m
        = \frac{3v\,q^2}{2(1-q)}
        \le 3vq^2
        = \frac{3r^2}{v},
    \end{equation*}
    using $1-q\ge\tfrac12$ in the last inequality. This proves (a).

    For (b), suppose $r^2\le v/12$. Since $v\ge1$, we have
    $r\le\sqrt{v/12}\le v/2$ (the latter being equivalent to $v\ge\tfrac13$), so
    part (a) applies and gives $|F(s)|\le 3r^2/v\le\tfrac14$. Setting
    $\theta:=e^{F(s)}-1$ and using $|e^z-1|\le|z|e^{|z|}$,
    \begin{equation*}
        |\theta| \le |F(s)|\,e^{|F(s)|}
        \le e^{1/4}\cdot\frac{3r^2}{v}
        \le \frac{4r^2}{v},
    \end{equation*}
    which is \eqref{eq:gamma_local_mult}; note that $|\theta|\le\tfrac13$.
    Finally, $\bigl||1+\theta|-1\bigr|\le|\theta|$ by the triangle inequality,
    and since $|\theta|\le1$ the point $1+\theta$ lies in the closed disc of
    radius $|\theta|$ centered at $1$, whence
    $|\arg(1+\theta)|\le\arcsin|\theta|\le\tfrac{\pi}{2}|\theta|$, the last step
    by convexity of $\arcsin$ on $[0,1]$.
\end{proof}

Specializing $w$ to the zero $\rho_j$ yields the bound we need for
$\tilde c_\rho$.

\begin{corollary}[Local coefficient bound]\label{cor:local_coeff}
    Fix any $\eps_1>0$,
    and fix a non-trivial zero with $\gamma_j\ge48\eps_1^{-2}$.
    Then for all non-trivial zeros $\rho$ with $|\rho-\rho_j|\le2\eps_1^{-1}$,
    \begin{equation}\label{eq:local_coeff}
        \frac{\tilde c_\rho}{\Gamma(\rho_j)} = 1 + \theta_\rho,
        \qquad \text{where $\theta_\rho\in\mathbb C$ satisfies }
        |\theta_\rho| \le \frac{C}{\eps_1^{2}\gamma_j}.
    \end{equation}
    In particular, \eqref{eq:local_coeff} controls both the modulus and the
    argument of $\tilde c_\rho/\Gamma(\rho_j)$, and $C=16$.
\end{corollary}

\begin{proof}
    We prove this by applying Lemma~\ref{lem:gamma_local}(b) with $w=\rho_j$ and $s=\rho$. 
    The assumptions hold by $u=\beta_j\in(0,1)$, $v=\gamma_j>0$, and $|s-w|^2\le4\eps_1^{-2}\le\gamma_j/12$.
    By \eqref{eq:tildec_def}, the left-hand side of \eqref{eq:gamma_local_mult} is $\tilde c_\rho/\Gamma(\rho_j)$,
    and \eqref{eq:gamma_local_mult} gives
    $|\theta_\rho|\le4\,(2\eps_1^{-1})^{2}/\gamma_j=16\eps_1^{-2}\gamma_j^{-1}$.
\end{proof}

We now look at the global behavior of $\tilde c_{\rho}$ away from $\rho_j$.

\begin{lemma}[Global coefficient bound]\label{lem:global_coeff}
    Let $\gamma_j$ be sufficiently large and $\rho=\beta+i\gamma$ be any nontrivial zero of $\zeta$. We write $b_{\rho}=\gamma-\gamma_j$. Then there exists an absolute constant $C' > 0$ such that the following bounds hold.
    \begin{enumerate}
        \item[\textup{(a)}] For $\gamma > 0$,
        \begin{equation}\label{eq:coeff_exact}
            \frac{|\tilde c_\rho|}{|\Gamma(\rho_j)|} \asymp \left(\frac{\gamma}{\gamma_j}\right)^{\beta-1/2}\exp\left(O\left(\frac{|b_\rho|}{\gamma_j}\right)\right).
        \end{equation}
        \item[\textup{(b)}] For $\gamma \ge \gamma_j$,
        \begin{equation}\label{eq:global_coeff_above}
            |\tilde c_\rho| \ll |\Gamma(\rho_j)|\left(\frac{\gamma}{\gamma_j}\right)^{1/2}\exp\left(\frac{C'|b_\rho|}{\gamma_j}\right).
        \end{equation}
        \item[\textup{(c)}] For $0 < \gamma < \gamma_j$,
        \begin{equation}\label{eq:global_coeff_below}
            |\tilde c_\rho| \ll |\Gamma(\rho_j)|\left(\frac{\gamma_j}{\gamma}\right)^{1/2}\exp\left(\frac{C'|b_\rho|}{\gamma_j}\right).
        \end{equation}
        \item[\textup{(d)}] For $\gamma<0$,
        \begin{equation}\label{eq:global_coeff_neg}
            |\tilde c_\rho| \ll |\Gamma(\rho_j)|\,\gamma_j^{1/2}e^{-\pi|\gamma|/2}.
        \end{equation}
    \end{enumerate}
\end{lemma}

\begin{proof}
    We first prove (a). By the Stirling formula, $|\Gamma(\rho)| \asymp \gamma^{\beta-1/2}e^{-\pi\gamma/2}$ and $|\Gamma(\rho_j)| \asymp \gamma_j^{\beta_j-1/2}e^{-\pi\gamma_j/2}$. Since $H_j(\rho - \rho_j) = (\sigma_j + i\nu_j)(a_\rho + ib_\rho)$, we have
    \begin{equation*}
        \Re\bigl(H_j(\rho - \rho_j)\bigr) = \sigma_j a_\rho - \nu_j b_\rho = a_\rho\log\gamma_j - \tfrac{\pi}{2}b_\rho + O\bigl((|a_\rho| + |b_\rho|)/\gamma_j\bigr)
    \end{equation*}
    by \eqref{eq:digamma_asym}. Therefore
    \begin{equation*}
        |\tilde c_\rho| = |\Gamma(\rho)|e^{-\Re(H_j(\rho-\rho_j))} \asymp \gamma^{\beta-1/2}e^{-\pi\gamma/2}\cdot e^{-a_\rho\log\gamma_j + (\pi/2)b_\rho}\cdot e^{O(|b_\rho|/\gamma_j)},
    \end{equation*}
    since the term $O(|a_\rho|/\gamma_j) = O(1/\gamma_j)$ is bounded by a constant (since $|a_\rho| \le 1$).
    Taking the ratio with $|\Gamma(\rho_j)| \asymp \gamma_j^{\beta_j - 1/2}e^{-\pi\gamma_j/2}$, we see that the exponential factors collapse into
    \begin{equation*}
        e^{-\pi\gamma/2}\cdot e^{(\pi/2)b_\rho}\cdot e^{\pi\gamma_j/2} = e^{-\pi\gamma/2 + (\pi/2)(\gamma-\gamma_j) + \pi\gamma_j/2} = 1
    \end{equation*}
    (since $b_\rho = \gamma - \gamma_j$), and the powers of $\gamma_j$ become $\gamma_j^{-(\beta_j-1/2)}\cdot\gamma_j^{-a_\rho} = \gamma_j^{-(\beta-1/2)}$ (using $\beta = \beta_j + a_\rho$). This proves \eqref{eq:coeff_exact}.

    Parts (b) and (c) follow from (a), since $0 < \beta < 1$, the exponent $\beta - 1/2 \in (-1/2, 1/2)$. For $\gamma \ge \gamma_j$ we have $\gamma/\gamma_j \ge 1$, and the worst case is $\beta - 1/2 = 1/2$; for $0 < \gamma < \gamma_j$ we have $\gamma/\gamma_j < 1$, and $(\gamma/\gamma_j)^{\beta-1/2} \le (\gamma_j/\gamma)^{1/2}$.

    For (d), we use the same computation as in (a), with $|\Gamma(\rho)|\asymp|\gamma|^{\beta-1/2}e^{-\pi|\gamma|/2}$ and $b_\rho=-|\gamma|-\gamma_j$, so
    \begin{equation*}
        \frac{|\tilde c_\rho|}{|\Gamma(\rho_j)|}
        \asymp\Bigl(\frac{|\gamma|}{\gamma_j}\Bigr)^{\beta-1/2}
        e^{-\pi|\gamma|}\,e^{O(1)+O(|\gamma|/\gamma_j)}.
    \end{equation*}
    Since $\beta-\tfrac12\in(-\tfrac12,\tfrac12)$ and $|\gamma|\ge10$, the first factor is at most $\gamma_j^{1/2}+|\gamma|^{1/2}$, and
    $|\gamma|^{1/2}e^{O(|\gamma|/\gamma_j)}e^{-\pi|\gamma|}\le e^{-\pi|\gamma|/2}$
    for sufficiently large $\gamma_j$.
    This proves \eqref{eq:global_coeff_neg}.
\end{proof}

\subsection{The main transform and a reduction proposition}
\label{subsec:kernel_detail}

We conclude some variations of Pintz~\cite[Section 5 to Section 6]{Pintz2017}.
\begin{lemma}\label{lem:mellin_gauss}
Fix $\rho_j = \beta_j + i\gamma_j$ with $\gamma_j > 0$ large, let $\eps_1 > 0$, $\alpha_j := \eta(\gamma_j)$, $L_j := \pi/(2|\eta'(\gamma_j)|)$, $A_j := \alpha_j L_j$, and let $\mu$ range over $\mathcal I_j := [L_j - 6\eps_1 A_j, L_j - 5\eps_1 A_j]$. Set $k := 5\eps_1^2\alpha_j\mu$ and $\mu_c := \mu - H_j$. For each zero $\rho$ of $\zeta$, let
\begin{equation}\label{eq:lambda_def}
    \lambda_\rho := 5\eps_1^2\alpha_j(\rho-\rho_j)^2 + (\rho-\rho_j),
\end{equation}
and define
\begin{equation}\label{eq:Uj_def}
    U_j(\mu) := -\sum_\rho\tilde c_\rho e^{\lambda_\rho\mu},
\end{equation}
where $\tilde c_\rho$ is as in \eqref{eq:tildec_def}. Using $\mu_c = \mu - H_j$,
\begin{equation*}
    U_j(\mu) = -\sum_\rho\Gamma(\rho)\exp\bigl(k(\rho-\rho_j)^2 + \mu_c(\rho-\rho_j)\bigr),
\end{equation*}
which treats $U_j(\mu)$ as an exponential-weighted (or Gaussian) sum over the zeros centered at $\rho_j$.

    Then, for $\mu \in \mathcal I_j$ and $j$ sufficiently large,
    \begin{equation}\label{eq:U_integral}
        U_j(\mu) = \frac{1}{\sqrt{4\pi k}}\int_0^\infty\Delta(e^t)e^{-\rho_j t}\exp\left(-\frac{(t-\mu_c)^2}{4k}\right)dt + O(e^{-cL_j})
    \end{equation}
    for some constant $c \in (0,1)$ (independent of $j$ and $\mu$).
\end{lemma}

\begin{proof}
    The explicit formula gives $\Delta(e^t) = -\sum_\rho\Gamma(\rho)e^{\rho t} + O(1)$ for $t \ge t_0$. Consider
    \begin{equation*}
        I(\mu) := \frac{1}{\sqrt{4\pi k}}\int_{-\infty}^\infty\Delta(e^t)e^{-\rho_j t}e^{-(t-\mu_c)^2/(4k)}\,dt.
    \end{equation*}

    Writing $s = \rho - \rho_j$ and using
    \begin{equation*}
        (t-\mu_c)^2 = (t-\Re\mu_c)^2 - (\Im\mu_c)^2 - 2i(t-\Re\mu_c)\Im\mu_c,
    \end{equation*}
    the absolute value of the integrand is
    \begin{equation*}
        |e^{st}|\cdot|e^{-(t-\mu_c)^2/(4k)}| = e^{t\Re s}\cdot e^{-((t-\Re\mu_c)^2 - (\Im\mu_c)^2)/(4k)}.
    \end{equation*}
    Integrating in $t$ by completing the square,
    \begin{equation*}
        \frac{1}{\sqrt{4\pi k}}\int_{-\infty}^\infty|e^{st}|\cdot\bigl|e^{-(t-\mu_c)^2/(4k)}\bigr|\,dt = \exp\bigl(k(\Re s)^2 + \Re\mu_c\cdot\Re s + (\Im\mu_c)^2/(4k)\bigr).
    \end{equation*}
    Now $|\Re s| = |\beta - \beta_j| \le 1$, $|\Re\mu_c| = |\mu - \sigma_j| \le \mu$, and $(\Im\mu_c)^2/(4k) = \nu_j^2/(4k) = O(1)$. Hence, for fixed $j$ and $\mu \in \mathcal I_j$, the right-hand side is bounded uniformly in $\rho$ by $\exp(O(k+\mu))$ (the implicit constant and the threshold in $j$ being absolute). Therefore
    \begin{equation*}
        \sum_\rho|\Gamma(\rho)|\cdot\exp(O(k+\mu)) = \exp(O(k+\mu))\sum_\rho|\Gamma(\rho)|.
    \end{equation*}
    By the Stirling formula, $|\Gamma(\rho)| \ll \gamma^{1/2}e^{-\pi\gamma/2}$, and the Riemann--von Mangoldt formula $N(T) \ll T\log T$ gives $\sum_\rho|\Gamma(\rho)| \ll \int_1^\infty t^{1/2}e^{-\pi t/2}\log t\,dt < \infty$.
    Thus we can then change the order of the integral and the sum.
    The identity (see: Pintz \cite[Equation (5.7)]{Pintz2017})
    \begin{equation*}
        \frac{1}{\sqrt{4\pi k}}\int_{-\infty}^\infty e^{st}e^{-(t-\mu_c)^2/(4k)}\,dt = e^{ks^2+\mu_c s}
    \end{equation*}
    with $s = \rho-\rho_j$ gives the zero sum $U_j(\mu)$.

    We next bound the remainder from the $O(1)$ terms. The constant (bounded) part of $\Delta(e^t)$ contributes
    \begin{equation*}
        \frac{O(1)}{\sqrt{4\pi k}}\int_0^\infty e^{-\rho_j t}e^{-(t-\mu_c)^2/(4k)}\,dt.
    \end{equation*}
    Extending the integral to $(-\infty,\infty)$ (the error will be handled in the coming paragraph), the Gaussian identity gives the integral $e^{k\rho_j^2 - \rho_j(\mu-\sigma_j) + i\nu_j\rho_j}$ up to the $\sqrt{4\pi k}$ factor, this can be computed by completing the square in $t$. Its absolute value is $\exp(\Re(k\rho_j^2 - \rho_j(\mu-\sigma_j) + i\nu_j\rho_j))$, which expands to
    \begin{equation*}
        k\beta_j^2 - k\gamma_j^2 - \beta_j(\mu-\sigma_j) - \nu_j\gamma_j \le k\beta_j^2 - \beta_j(\mu-\sigma_j) + O(1),
    \end{equation*}
    after dropping the two negative terms $-k\gamma_j^2$ and $-\nu_j\gamma_j = -(\pi/2)\gamma_j + O(1)$. The remaining terms are controlled as follows: $k\beta_j^2 = 5\eps_1^2\alpha_j\mu\beta_j^2 = o(\mu)$ since $\alpha_j\to 0$, and $\beta_j\sigma_j = \beta_j\log\gamma_j = o(\mu)$ since $\gamma_j = o(A_j) = o(\mu)$. Thus the exponent is $\le -\beta_j\mu + o(\mu) \le -\mu/2$ for $\beta_j \ge 3/4$ (which holds for $j$ large by \eqref{eq:violate}). Since $\mu \asymp L_j$, the total remainder is $O(e^{-\mu/2}) = O(e^{-L_j/4})$.

    We finally handle the error of truncation to $[0,\infty)$. For $t \le 0$, $|\Delta(e^t)| = |\Delta(u)|$ with $u = e^t \in (0,1]$, so $|\Delta(e^t)| = O(1)$. The modulus of the integrand is bounded by
    \begin{equation*}
        O(1)\cdot e^{-\beta_j t}\cdot\exp\left(-\frac{(t-\Re\mu_c)^2 - (\Im\mu_c)^2}{4k}\right), \qquad \Re\mu_c = \mu - \sigma_j,\ \Im\mu_c = -\nu_j.
    \end{equation*}
    Instead of dropping $e^{-\beta_j t}$ on $(-\infty,0]$, we complete the square in $t$ by combining the $e^{-\beta_j t}$ factor into the exponential part,
    \begin{equation*}
        \begin{aligned}
            -\beta_j t - \frac{(t-\Re\mu_c)^2}{4k} &= -\frac{(t-t_0)^2}{4k} + k\beta_j^2 - \beta_j\Re\mu_c, \\
            t_0 &:= \Re\mu_c - 2k\beta_j = \mu - \sigma_j - 2k\beta_j.
        \end{aligned}
    \end{equation*}
    Since $\mu \ge L_j - 6\eps_1 A_j \sim L_j$, $\sigma_j = \log\gamma_j = o(L_j)$, and $2k\beta_j \le 2k = O(A_j) = o(L_j)$, we have $t_0 \ge L_j/2$ for $j$ large, hence for $t \le 0 \le t_0$, $(t-t_0)^2 \ge t_0^2 \ge L_j^2/4$, so $(t-t_0)^2/(4k) \ge L_j^2/(16k)$. The $(-\infty,0]$ contribution is therefore bounded by
    \begin{equation*}
        \begin{split}
            &O(1)\cdot\frac{1}{\sqrt{4\pi k}}e^{k\beta_j^2 - \beta_j\Re\mu_c + (\Im\mu_c)^2/(4k)}\cdot e^{-L_j^2/(16k)}\int_{-\infty}^0\exp\left(-\frac{(t-t_0)^2 - L_j^2/4}{4k}\right)dt \\
            &\le O(1)\cdot e^{k\beta_j^2 - L_j^2/(16k) + O(1)},
        \end{split}
    \end{equation*}
    using $e^{-\beta_j\Re\mu_c} \le 1$ (since $\Re\mu_c > 0$) and $(\Im\mu_c)^2/(4k) = \nu_j^2/(4k) = O(1/k) = o(1)$. Here the residual Gaussian integral (centered at $t_0$ with variance $2k$, restricted to $(-\infty,0]$) is at most the full Gaussian mass $\sqrt{4\pi k}$, which cancels the $1/\sqrt{4\pi k}$ factor.

    It remains to bound the exponent $k\beta_j^2 - L_j^2/(16k)$. Since $k = o(L_j)$, for $j$ large $k \le L_j/32$, giving $L_j^2/(16k) \ge 2L_j$; and $k\beta_j^2 \le k \le L_j/32 \le L_j$. Therefore
    \begin{equation*}
        k\beta_j^2 - L_j^2/(16k) \le L_j - 2L_j = -L_j,
    \end{equation*}
    so the $(-\infty,0]$ part is $O(e^{-L_j + O(1)}) = O(e^{-cL_j})$ for any $c\in(0,1)$ and $j$ large.
\end{proof}

Finally, we prove the following proposition, which together with Theorem~\ref{thm:main1} yields Theorem~\ref{thm:main2}.

\begin{proposition}\label{prop:Gk}
    Let $k \ge 2$ and let $\varphi:[1,\infty)\to\mathbb R$ satisfy $\varphi(x) \to \infty$ and $xe^{-\varphi(x)} \to \infty$ as $x $ tends to infinity. Then\\
    \textup{(a)} $\Delta(x) \ll xe^{-\varphi(x)}$ implies $G_k(x)-x^k \ll_k x^k e^{-\varphi(x)}$;\\
    \textup{(b)} $G_k(x)-x^k \ll_k x^k e^{-\varphi(x)}$ implies $\Delta(x) \ll_k xe^{-\varphi(x)}$.
\end{proposition}

\begin{proof}
    Since $G_k(x) = \Psi(x)^k$ and $\Psi(x) = x + \widetilde\Delta(x)$ with $\widetilde\Delta(x) := \Delta(x) - \frac{1}{2} + O(1/x)$, the hypothesis $xe^{-\varphi(x)} \to \infty$ makes the $O(1)$ difference between $\Delta$ and $\widetilde\Delta$ negligible compared with $xe^{-\varphi(x)}$:
    \begin{equation*}
        |\Delta(x)| \ll xe^{-\varphi(x)} \quad\Longleftrightarrow\quad |\widetilde\Delta(x)| \ll xe^{-\varphi(x)}.
    \end{equation*}

    For the forward direction, assume $|\widetilde\Delta(x)| \le Cxe^{-\varphi(x)}$, so $|\widetilde\Delta(x)/x| \le Ce^{-\varphi(x)} \to 0$. Expanding by the binomial theorem,
    \begin{equation*}
        G_k(x) - x^k = (x + \widetilde\Delta(x))^k - x^k = \sum_{j=1}^k\binom{k}{j}x^{k-j}\widetilde\Delta(x)^j.
    \end{equation*}
    Taking the $j=1$ term $kx^{k-1}\widetilde\Delta(x)$ out of the sum,
    \begin{equation*}
        G_k(x) - x^k = kx^{k-1}\widetilde\Delta(x)\left(1 + \sum_{j=2}^k\frac{1}{k}\binom{k}{j}\left(\frac{\widetilde\Delta(x)}{x}\right)^{j-1}\right),
    \end{equation*}
    using $\binom{k}{j}x^{k-j}\widetilde\Delta^j/(kx^{k-1}\widetilde\Delta) = \binom{k}{j}(\widetilde\Delta/x)^{j-1}/k$. Since $\widetilde\Delta(x)/x \to 0$, the inner sum is bounded by $C_k|\widetilde\Delta(x)/x| \le C_ke^{-\varphi(x)}$ for $x$ large. Hence
    \begin{equation*}
        |G_k(x) - x^k| \le kx^{k-1}\cdot Cxe^{-\varphi(x)}\cdot(1 + C_ke^{-\varphi(x)}) \le C'_k x^ke^{-\varphi(x)}.
    \end{equation*}

    In the reverse direction, assume $|G_k(x)-x^k| \le C_kx^ke^{-\varphi(x)}$. Since $G_k(x) = \Psi(x)^k = x^k(1+O(e^{-\varphi(x)}))$ with $\varphi \to \infty$, $\Psi(x) > 0$ for large $x$ and
    \begin{equation*}
        \frac{\Psi(x)}{x} = (1+O(e^{-\varphi(x)}))^{1/k} = 1 + O_k(e^{-\varphi(x)}),
    \end{equation*}
    using $(1+u)^{1/k} = 1+O(u)$ for $|u| \le 1/2$. Hence $\widetilde\Delta(x) = \Psi(x) - x = O_k(xe^{-\varphi(x)})$.
\end{proof}

\section{Proofs of the Main Theorems}
\label{sec:proofs}

\subsection{Proof of Theorem~\ref{thm:main1}: Upper Bound}

Define
\begin{align*}
    \omega(x) &:= \min_{\substack{\rho=\beta+i\gamma\\\gamma>0}}\bigl((1-\beta)\log x - \log|\Gamma(\rho)|\bigr), \\
    S(x) &:= \max_{u\le x}|\Delta(u)|, \\
    W(x) &:= \sum_{0<\gamma\le x}\frac{|\Gamma(\rho+1)|x^\beta}{\gamma}.
\end{align*}

\begin{proof}
    Assume $\zeta(\sigma+it) \neq 0$ for $\sigma > 1-\eta(|t|)$, $|t| \ge 2$. Since the first non-trivial zero of $\zeta$ satisfies $\gamma_1 > 10$, every non-trivial zero $\rho = \beta+i\gamma$ with $\gamma > 0$ has $\gamma \ge 10$, and in particular $1-\beta \ge \eta(\gamma)$.

    We begin with the following bound. By \eqref{eq:stirling_master}, uniformly for $0 \le \beta \le 1$,
    \begin{align*}
        \log|\Gamma(\beta+i\gamma)|
        &= (\beta-\tfrac{1}{2})\log\gamma - \tfrac{\pi\gamma}{2} + \tfrac{1}{2}\log(2\pi) + O(\gamma^{-1}) \\
        &= \log|\Gamma(i\gamma)| + \beta\log\gamma + O(\gamma^{-1}).
    \end{align*}
    Hence for any non-trivial zero $\rho$ (so $\gamma \ge 10$, $1-\beta \ge \eta(\gamma)$, and $\beta \le 1$),
    \begin{align}
        (1-\beta)\log x - \log|\Gamma(\rho)|
        &= (1-\beta)\log x - \log|\Gamma(i\gamma)| - \beta\log\gamma + O(\gamma^{-1}) \notag \\
        &\ge \eta(\gamma)\log x - \log|\Gamma(i\gamma)| - \log\gamma + O(1) \notag \\
        &\ge \omega_\eta(x) - \log\gamma + O(1). \label{eq:pointwise_omega}
    \end{align}

    We next verify that the infimum $\omega(x)$ is attainable. For fixed $x$, $(1-\beta)\log x - \log|\Gamma(\rho)| \ge -\log|\Gamma(\rho)| \sim \tfrac{\pi}{2}\gamma$ as $\gamma \to \infty$, so only finitely many zeros can give a value below any threshold; the infimum is therefore attained at some $\rho_0 = \beta_0 + i\gamma_0$, with $\gamma_0 \ge 10$.

    We now bound $\gamma_0$. Since $\beta_0 \le 1$, $(1-\beta_0)\log x \ge 0$, hence
    \begin{equation*}
        \omega(x) = (1-\beta_0)\log x - \log|\Gamma(\rho_0)| \ge -\log|\Gamma(\rho_0)| = \tfrac{\pi}{2}\gamma_0 + O(\log\gamma_0)
    \end{equation*}
    by the Stirling formula. For $\gamma_0$ large enough that the $O(\log\gamma_0)$ term is absorbed, $\omega(x) \ge \tfrac{\pi}{3}\gamma_0$, i.e.\ $\gamma_0 \le 3\omega(x)/\pi$, and therefore
    \begin{equation}\label{eq:gamma0_bound}
        \log\gamma_0 \le \log\omega(x) + O(1).
    \end{equation}

    Combining the above, we apply \eqref{eq:pointwise_omega} to $\rho = \rho_0$:
    \begin{equation*}
        \begin{aligned}
            \omega(x) &= (1-\beta_0)\log x - \log|\Gamma(\rho_0)| \\
            &\ge \omega_\eta(x) - \log\gamma_0 + O(1) \\
            &\ge \omega_\eta(x) - \log\omega(x) + O(1),
        \end{aligned}
    \end{equation*}
    where the second inequality uses \eqref{eq:gamma0_bound}.
    We first note that $\omega_\eta(x)\to\infty$: for any fixed $T\ge2$, splitting the infimum in \eqref{eq:omega_eta} at $t=T$ and using that $\eta$ is decreasing and $g$ is increasing gives
    \begin{equation*}
        \omega_\eta(x)\ \ge\ \min\bigl(\eta(T)\log x+g(2),\ g(T)\bigr);
    \end{equation*}
    letting $x\to\infty$ and then $T\to\infty$ proves the claim.
    Now, if $\omega(x)\ge\omega_\eta(x)$, then \eqref{eq:omega_compare} holds trivially.
    Otherwise $\omega(x)<\omega_\eta(x)$, so $\log\omega(x)\le\log\omega_\eta(x)$
    (notice that $\omega(x)\ge-\log|\Gamma(\rho_0)|\ge\pi\gamma_1/2-O(\log\gamma_1)>1$), and the previous display yields
    \begin{equation*}
        \omega(x)\ \ge\ \omega_\eta(x)-\log\omega_\eta(x)+O(1)
        \ =\ (1-o(1))\,\omega_\eta(x),
    \end{equation*}
    since $\log\omega_\eta(x)=o(\omega_\eta(x))$ as $\omega_\eta(x)\to\infty$.
    Hence
    \begin{equation}\label{eq:omega_compare}
        \omega(x) \ge \omega_\eta(x)(1-o(1)).
    \end{equation}

    We finally pass the results from $\omega$ to $\Delta$. By the explicit formula \eqref{explicit_formula}, for $u \le x$,
    \begin{equation*}
        |\Delta(u)| \le \sum_\rho|\Gamma(\rho)|u^\beta + O(1) \le \sum_\rho|\Gamma(\rho)|x^\beta + O(1).
    \end{equation*}
    Using $|\Gamma(\rho)| = |\Gamma(\rho+1)|/|\rho| \le |\Gamma(\rho+1)|/\gamma$ and splitting the sum at $\gamma = x$, zeros with $\gamma > x$ contribute $O(e^{-x})$ by the exponential decreasing of $|\Gamma(\rho)|$. Hence $S(x) \ll W(x) + O(1)$, and $W(x) \to \infty$.

    The upper bound $W(x) \ll xe^{-(1-\eps)\omega(x)}$ is proved in Section~\ref{sec:upper_bound}. Combining these we have,
    \begin{equation*}
        |\Delta(x)| \le S(x) \ll W(x) \ll xe^{-(1-\eps)\omega(x)} \ll xe^{-(1-\eps - o(1))\omega_\eta(x)}.
    \end{equation*}
    Since $\eps > 0$ is arbitrary, replacing $\eps$ by $\eps/2$ and absorbing the $o(1)$ for $x$ sufficiently large yields $\Delta(x) \ll xe^{-(1-\eps)\omega_\eta(x)}$ for any fixed $\eps \in (0,1)$.
\end{proof}

\subsection{Proof of Theorem~\ref{thm:main1}: Converse}
\label{subsec:converse}

The converse is proved by contradiction using a Tur\'an power-sum argument, using an improved method of Pintz~\cite{Pintz2017}, applied directly to the violating zeros.

\begin{proof}
    Assume $\Delta(x) \ll xe^{-(1+\eps)\omega_\eta(x)}$. Fix $\delta > 0$ and suppose, by contradiction, that the claimed ZFR with $(1-\delta)\eta$ fails for arbitrarily large $t$. Then there exists a sequence of zeros $\{\rho_j = \beta_j + i\gamma_j\}_j$ with $\gamma_j \to \infty$ and
    \begin{equation}\label{eq:violate}
        1-\beta_j < (1-\delta)\eta(\gamma_j).
    \end{equation}
    Since the final contradiction is derived for each fixed $\delta>0$, this yields the conclusion of the theorem.

    With $\delta>0$ fixed as above and $\eps>0$ as in the hypothesis, we fix the parameter
    \begin{equation*}
        \eps_1 := \frac{\min(\eps,\delta)}{48}
    \end{equation*}
    and set $\alpha_j := \eta(\gamma_j) \to 0$. Note that $\eps_1$ depends on both $\eps$ (fixed once and for all by the hypothesis) and $\delta$ (fixed at the beginning of the contradiction argument). In particular, $\eps_1 \le \eps/48$ and $\eps_1 \le \delta/48$. Define
    \begin{equation*}
        L_j := \frac{\pi}{2|\eta'(\gamma_j)|}, \qquad A_j := \alpha_j L_j, \qquad B_j := g(\gamma_j) = -\log|\Gamma(i\gamma_j)|.
    \end{equation*}
    By (H3), $A_j/\gamma_j \to \infty$, in particular $A_j \to \infty$ and $B_j = o(A_j)$ since $B_j \sim \frac{\pi}{2}\gamma_j$. Let $\mu$ be in the interval
    \begin{equation}\label{eq:I_j}
        \mathcal{I}_j := [L_j - 6\eps_1 A_j, L_j - 5\eps_1 A_j],
    \end{equation}
    and define
    \begin{equation}\label{eq:kM}
        k = k(\mu) := 5\eps_1^2\alpha_j\mu, \qquad M = M(\mu) := k/\eps_1 = 5\eps_1\alpha_j\mu.
    \end{equation}
    For $\mu \in \mathcal{I}_j$ we have $k \asymp \eps_1^2 A_j \to \infty$, $M \asymp \eps_1 A_j$, and $\mu + M = \mu(1+5\eps_1\alpha_j) \le L_j$, and so $e^{\mu+M} \le x_j := e^{L_j}$. This is because $\mu \le L_j - 5\eps_1 A_j$ and $5\eps_1\alpha_j\mu \le 5\eps_1 A_j$.

    We next set up the kernel and its integral representation. With $H_j$ as in \eqref{eq:digamma_asym} and $\mu_c := \mu - H_j$, let $U_j(\mu)$ and $\lambda_\rho$ be defined by \eqref{eq:Uj_def} and \eqref{eq:lambda_def}, so
    \begin{equation}\label{eq:U_power_sum}
        U_j(\mu) = -\sum_\rho\tilde c_\rho e^{\lambda_\rho\mu},
    \end{equation}
    where $\tilde c_\rho$ is as in \eqref{eq:tildec_def}, note $\lambda_{\rho_j} = 0$ and $\tilde c_{\rho_j} = \Gamma(\rho_j)$. By Lemma~\ref{lem:mellin_gauss}, for $\mu \in \mathcal{I}_j$,
    \begin{equation*}
        U_j(\mu) = \frac{1}{\sqrt{4\pi k}}\int_0^\infty \Delta(e^t)e^{-\rho_j t}\exp\left(-\frac{(t-\mu_c)^2}{4k}\right)dt + O(e^{-cL_j}).
    \end{equation*}

    We now bound the integral from above. Split the integral into three parts, which we denote as follows. Set
    \begin{equation*}
        F_j(t;\mu) := \frac{1}{\sqrt{4\pi k}}\Delta(e^t)e^{-\rho_j t}
        \exp\left(-\frac{(t-\mu_c)^2}{4k}\right),
    \end{equation*}
    and define
    \begin{align*}
        U_{1,j}(\mu) &:= \int_0^{\mu-M} F_j(t;\mu)\,dt, \\
        U_{2,j}(\mu) &:= \int_{\mu-M}^{\mu+M} F_j(t;\mu)\,dt, \\
        U_{3,j}(\mu) &:= \int_{\mu+M}^{\infty} F_j(t;\mu)\,dt.
    \end{align*}
    Thus $U_j(\mu)=U_{1,j}(\mu)+U_{2,j}(\mu)+U_{3,j}(\mu)+O(e^{-cL_j})$.
    We suppress the argument $\mu$ below.

    We estimate the central integral $U_{2,j}$ first. On $[\mu-M, \mu+M]$ we have $|t - \mu| \le M$ and $\mu_c = \mu - \sigma_j - i\nu_j$, so $t - \mu_c = (t - \mu + \sigma_j) + i\nu_j$. Expanding the square,
    \begin{equation*}
        (t-\mu_c)^2 = (t-\mu+\sigma_j)^2 - \nu_j^2 + 2i\nu_j(t-\mu+\sigma_j),
    \end{equation*}
    whose real part is $(t-\mu+\sigma_j)^2 - \nu_j^2 \ge -\nu_j^2$. Therefore
    \begin{equation*}
        \left|\exp\left(-\frac{(t-\mu_c)^2}{4k}\right)\right| = \exp\left(-\frac{(t-\mu+\sigma_j)^2 - \nu_j^2}{4k}\right) \le \exp\left(\frac{\nu_j^2}{4k}\right) = 1 + o(1)
    \end{equation*}
    since $\nu_j = O(1)$ and $k \to \infty$. Thus
    \begin{equation*}
        |U_{2,j}| \ll \frac{1}{\sqrt{k}}\int_{\mu-M}^{\mu+M}|\Delta(e^t)|e^{-\beta_j t}\,dt.
    \end{equation*}
    Substituting $u = e^t$ (so $du/u = dt$ and $e^{-\beta_j t} = u^{-\beta_j}$),
    \begin{equation*}
        |U_{2,j}| \ll \frac{1}{\sqrt{k}}\int_{e^{\mu-M}}^{e^{\mu+M}}\frac{|\Delta(u)|}{u^{1+\beta_j}}\,du \le \frac{e^{-\beta_j(\mu-M)}}{\sqrt{k}}\int_{e^{\mu-M}}^{e^{\mu+M}}\frac{|\Delta(u)|}{u}\,du,
    \end{equation*}
    where we used $u^{-\beta_j} \le e^{-\beta_j(\mu - M)}$ on the interval of integration (the lower endpoint is where $u$ is smallest). For the remaining integral, on the interval $[e^{\mu-M}, e^{\mu+M}] \subseteq [0, x_j]$ (since $e^{\mu+M} \le x_j$) we have $u^{-1} \le e^{-(\mu - M)}$, giving
    \begin{equation*}
        \int_{e^{\mu-M}}^{e^{\mu+M}}\frac{|\Delta(u)|}{u}\,du \le e^{-(\mu-M)}\int_0^{x_j}|\Delta(u)|\,du = e^{-(\mu-M)}x_j D(x_j),
    \end{equation*}
    where
    \begin{equation}\label{eq:D_def}
        D(x) := \frac{1}{x}\int_0^x|\Delta(u)|\,du.
    \end{equation}
    Combining the two estimates and writing $x_j = e^{L_j}$,
    \begin{equation*}
        |U_{2,j}| \ll \frac{1}{\sqrt{k}}e^{-\beta_j(\mu-M)}\cdot e^{-(\mu-M)}\cdot e^{L_j}D(x_j) = \frac{x_j^{-\beta_j}}{\sqrt{k}}e^{(1+\beta_j)(L_j-\mu+M)}D(x_j),
    \end{equation*}
    using $e^{-(1+\beta_j)(\mu-M)}\cdot e^{L_j} = e^{-\beta_j L_j}\cdot e^{(1+\beta_j)(L_j - \mu + M)} = x_j^{-\beta_j}e^{(1+\beta_j)(L_j-\mu+M)}$. From \eqref{eq:I_j} and \eqref{eq:kM}, $L_j - \mu + M \le 6\eps_1 A_j + 5\eps_1 A_j = 11\eps_1 A_j$. Since $1+\beta_j \le 2$,
    \begin{equation}\label{eq:U2_bound}
        |U_{2,j}| \ll \frac{x_j^{-\beta_j}}{\sqrt{k}}e^{22\eps_1 A_j}D(x_j).
    \end{equation}

    We next estimate the integral $U_{3,j}$ over $t \ge \mu+M$. Using $|\Delta(e^t)| \ll e^t$, which holds unconditionally since $\Psi(x) \ll x$ by Chebyshev's bound and $I(x) = O(x)$, the integrand is $O(e^{h(t)}/\sqrt{k})$ where $h(t) := (1-\beta_j)t - (t-\mu+\sigma_j)^2/(4k)$. Its derivative is
    \begin{equation*}
        h'(t) = (1-\beta_j) - \frac{t-\mu+\sigma_j}{2k}.
    \end{equation*}
    At $t = \mu+M$, $h'(\mu+M) \le (1-\beta_j) - M/(2k) = (1-\beta_j) - 1/(2\eps_1)$. Since $(1-\beta_j) \le \alpha_j \to 0$ while $1/(2\eps_1)$ is a fixed positive constant (determined by $\eps$ and $\delta$), $h'(\mu+M) < 0$ for large $j$. Since $h''(t) = -1/(2k) < 0$, $h'$ is decreasing, so $h'(t) < 0$ for all $t \ge \mu+M$. This part is therefore dominated by
    \begin{equation*}
        |U_{3,j}| \ll \frac{1}{\sqrt{k}}e^{h(\mu+M)}.
    \end{equation*}
    We drop the term $2M\sigma_j + \sigma_j^2$ and note that $h(\mu+M) \le (1-\beta_j)(\mu+M) - M^2/(4k)$. Substituting $M = k/\eps_1$ and $k = 5\eps_1^2\alpha_j\mu$ gives $M^2/(4k) = k/(4\eps_1^2) = (5/4)\alpha_j\mu$. Using \eqref{eq:violate} and $\mu+M \le L_j$, we have $(1-\beta_j)(\mu+M) \le (1-\delta)\alpha_j L_j = (1-\delta)A_j$. For the lower bound on $(5/4)\alpha_j\mu$, since $\mu \ge L_j - 6\eps_1 A_j$,
    \begin{equation*}
        \tfrac{5}{4}\alpha_j\mu \ge \tfrac{5}{4}\alpha_j(L_j - 6\eps_1 A_j) = \tfrac{5}{4}A_j(1-6\eps_1\alpha_j) \ge A_j
    \end{equation*}
    for $j$ large, since $(5/4)(1-6\eps_1\alpha_j) \to 5/4 > 1$ as $\alpha_j \to 0$. Therefore
    \begin{equation}\label{eq:tail_exponent}
        h(\mu+M) \le (1-\delta)A_j - A_j + o(A_j) = -\delta A_j + o(A_j) \le -(\delta/2)A_j
    \end{equation}
    for $j$ large. A symmetric argument gives $|U_{1,j}| \ll e^{-(\delta/2)A_j}/\sqrt{k}$.

    By the Stirling formula, $|\Gamma(\rho_j)| = e^{-B_j + O(\log\gamma_j)}$ with $B_j = o(A_j)$. Set $c_0 := \delta/3$. Then
    \begin{equation*}
        \begin{split}
            -\tfrac{\delta}{2}A_j &= \log|\Gamma(\rho_j)| + B_j + O(\log\gamma_j) - \tfrac{\delta}{2}A_j \\
            &= \log|\Gamma(\rho_j)| - \bigl(\tfrac{\delta}{2} - o(1)\bigr)A_j\\
            &\le \log|\Gamma(\rho_j)| - c_0 A_j
        \end{split}
    \end{equation*}
    for $j$ large (since $\delta/2 - o(1) \ge \delta/3$), so $e^{-(\delta/2)A_j} \le |\Gamma(\rho_j)|e^{-c_0 A_j}$. Combining these together we have
    \begin{equation}\label{eq:U_upper}
        |U_j(\mu)| \ll \frac{x_j^{-\beta_j}}{\sqrt{k}}e^{22\eps_1 A_j}D(x_j) + |\Gamma(\rho_j)|e^{-c_0 A_j}.
    \end{equation}

    We now analyze the coefficients. Write $a_\rho := \beta - \beta_j$ and $b_\rho := \gamma - \gamma_j$. The coefficients $\tilde c_\rho$ were analyzed in Section~\ref{subsec:coeff}:
    \begin{itemize}
        \item By Corollary~\ref{cor:local_coeff}, for every zero $\rho$ with $|\rho-\rho_j|\le2\eps_1^{-1}$, once $\gamma_j\ge48\eps_1^{-2}$ (which holds for $j$ large), we have $\tilde c_\rho=\Gamma(\rho_j)(1+\theta_\rho)$ with $|\theta_\rho|\le16\eps_1^{-2}\gamma_j^{-1}$. This controls both the modulus and the argument.
        \item Lemma~\ref{lem:global_coeff} provides three global bounds \eqref{eq:global_coeff_above}, \eqref{eq:global_coeff_below}, \eqref{eq:global_coeff_neg}, for the cases $\gamma \ge \gamma_j$, $0 < \gamma < \gamma_j$, and $\gamma < 0$ respectively. Since all zeros satisfy $\gamma \ge \gamma_1 > 10$, the factor $(\gamma_j/\gamma)^{1/2}$ appearing in \eqref{eq:global_coeff_below} is at most $(\gamma_j/10)^{1/2} = O(\gamma_j^{1/2})$.
    \end{itemize}

    We classify the zeros and bound the error terms. Since $\beta_j > 1-(1-\delta)\alpha_j > 1-\alpha_j$ from \eqref{eq:violate}, we have
    \begin{equation}\label{eq:a_bounds}
        a_\rho < (1-\delta)\alpha_j \quad(\text{for all }\rho), \qquad a_\rho > -2\alpha_j \quad (\text{if }\beta > 1-2\alpha_j).
    \end{equation}
    Recall $\Re(\lambda_\rho) = 5\eps_1^2\alpha_j(a_\rho^2 - b_\rho^2) + a_\rho$. We split the zeros into three sets.

    We first consider the set $E_1 := \{\rho : \beta \le 1-2\alpha_j\}$. Since $\beta_j > 1 - (1-\delta)\alpha_j > 1 - \alpha_j$, we have $a_\rho = \beta - \beta_j < (1-2\alpha_j) - (1-\alpha_j) = -\alpha_j$. We now split according to the size of $|b_\rho|$:
    \begin{itemize}
        \item If $|b_\rho| \ge |a_\rho|$, then $a_\rho^2 - b_\rho^2 \le 0$, so $\Re(\lambda_\rho) \le a_\rho \le -\alpha_j$.
        \item If $|b_\rho| < |a_\rho|$, then $a_\rho^2 - b_\rho^2 \le a_\rho^2 \le 1$ since $|a_\rho| \le 1$ from $\beta,\beta_j\in[0,1]$, so $5\eps_1^2\alpha_j(a_\rho^2-b_\rho^2) \le 5\eps_1^2\alpha_j \le \alpha_j/2$, by using $\eps_1 \le 1/48$, so $5\eps_1^2 \le 5/48^2 \le 1/2$. Hence $\Re(\lambda_\rho) \le \alpha_j/2 + a_\rho \le \alpha_j/2 - \alpha_j = -\alpha_j/2$.
    \end{itemize}
    In either case $\Re(\lambda_\rho) \le -\alpha_j/2$, so $|e^{\lambda_\rho\mu}| \le e^{-\alpha_j\mu/2} \le e^{-A_j/3}$, where the last step uses $\alpha_j\mu \ge \alpha_j(L_j - 6\eps_1 A_j) = A_j(1 - 6\eps_1\alpha_j) \ge 2A_j/3$ for $j$ large.

    The set $E_3 := \{\rho : \beta > 1-2\alpha_j,\ |b_\rho| > \eps_1^{-1}\}$ is handled similarly. By \eqref{eq:a_bounds}, $a_\rho < \alpha_j$ and $5\eps_1^2\alpha_j a_\rho^2 \le 20\eps_1^2\alpha_j^3 = o(\alpha_j)$, while $5\eps_1^2\alpha_j b_\rho^2 > 5\alpha_j$. Hence $\Re(\lambda_\rho) \le \alpha_j + o(\alpha_j) - 5\alpha_j \le -3\alpha_j$.

    The remaining set $E_2 := \{\rho : \beta > 1-2\alpha_j,\ |b_\rho| \le \eps_1^{-1}\}$ is local: its zeros satisfy $|\rho-\rho_j| \le 2\alpha_j + \eps_1^{-1} < 2\eps_1^{-1}$, so \eqref{eq:local_coeff} applies.

    For $i\in\{1,2,3\}$, we set
    \begin{equation*}
        E_i(\mu) := \sum_{\rho\in E_i}\tilde c_\rho e^{\lambda_\rho\mu},
    \end{equation*}
    so that $-U_j(\mu) = E_1(\mu) + E_2(\mu) + E_3(\mu)$ by \eqref{eq:U_power_sum}.

    \medskip We now bound $|E_1(\mu)| + |E_3(\mu)|$ by summing over integer intervals.
    We first bound the contribution of $E_3$ over $\gamma > 0$. For $\rho \in E_3$, $|a_\rho| < 2\alpha_j$ from \eqref{eq:a_bounds}, so $a_\rho^2 < 4\alpha_j^2$ and $5\eps_1^2\alpha_j a_\rho^2 < 20\eps_1^2\alpha_j^3 = o(\alpha_j)$, hence
    \begin{equation*}
        \Re(\lambda_\rho) \le a_\rho + 5\eps_1^2\alpha_j a_\rho^2 - 5\eps_1^2\alpha_j b_\rho^2 \le \alpha_j + o(\alpha_j) - 5\eps_1^2\alpha_j b_\rho^2.
    \end{equation*}
    Multiplying by $\mu \le L_j$ and using $\alpha_j\mu \le A_j$ together with $\alpha_j\mu \ge A_j(1-6\eps_1\alpha_j)$,
    \begin{equation*}
        \Re(\lambda_\rho)\mu \le A_j - 5\eps_1^2\alpha_j\mu b_\rho^2 + o(A_j) \le A_j - 4\eps_1^2 A_j b_\rho^2 + o(A_j),
    \end{equation*}
    where we used $5\eps_1^2\alpha_j\mu \ge 5\eps_1^2 A_j(1-6\eps_1\alpha_j) \ge 4\eps_1^2 A_j$ for $j$ large. Since $\beta > 1/2$ in $E_3$, for $\gamma \ge \gamma_j$, \eqref{eq:global_coeff_above} gives $|\tilde c_\rho| \ll |\Gamma(\rho_j)|(\gamma/\gamma_j)^{1/2}\exp(C'|b_\rho|/\gamma_j)$; for $\gamma < \gamma_j$, Lemma~\ref{lem:global_coeff}(a) with $\beta-1/2 > 0$ gives $|\tilde c_\rho| \ll |\Gamma(\rho_j)|\exp(C'|b_\rho|/\gamma_j)$. In both cases
    \begin{equation*}
        |\tilde c_\rho| \ll |\Gamma(\rho_j)|\max\bigl(1,(\gamma/\gamma_j)^{1/2}\bigr)\exp(C'|b_\rho|/\gamma_j).
    \end{equation*}
    Next we partition the zeros of $E_3$ into intervals of length $1$ in $|b_\rho|$: for each integer $m \ge \lfloor\eps_1^{-1}\rfloor$, consider the zeros with $m < |b_\rho| \le m+1$. By the Riemann--von Mangoldt formula, there are $O(\log(\gamma_j+m))$ such zeros. Define the weight function
    \begin{equation}\label{eq:weight_def}
        w(m) := \log(\gamma_j+m)\cdot\bigl((\gamma_j+m)/\gamma_j\bigr)^{1/2}, \qquad m \ge 0,
    \end{equation}
    which captures the combined contribution of the zero-counting factor $\log(\gamma_j+m)$ and the polynomial factor $(\gamma/\gamma_j)^{1/2}$ from \eqref{eq:global_coeff_above}. Then
    \begin{equation}\label{eq:E3_interval}
        |E_3(\mu)| \ll |\Gamma(\rho_j)|\sum_{m\ge\eps_1^{-1}}w(m)\exp\left(\frac{C'm}{\gamma_j} + A_j - 4\eps_1^2 A_j m^2\right).
    \end{equation}
    We verify that consecutive terms decreasing very fast. Let
    \begin{equation*}
        \phi(m) := \frac{C'm}{\gamma_j} + A_j - 4\eps_1^2 A_j m^2
    \end{equation*}
    denote the exponent at index $m$. The difference between indices $m+1$ and $m$ is
    \begin{equation}\label{eq:interval_diff}
        \phi(m+1) - \phi(m) = \frac{C'}{\gamma_j} - 4\eps_1^2 A_j(2m+1).
    \end{equation}
    For $m \ge \eps_1^{-1}$, $4\eps_1^2 A_j(2m+1) \ge 4\eps_1^2 A_j\cdot 2\eps_1^{-1} = 8\eps_1 A_j$, while $C'/\gamma_j = O(1/\gamma_j) = o(1)$. So $\phi(m+1) - \phi(m) \le -7\eps_1 A_j$ for $j$ large. The weight $w$ satisfies
    \begin{equation}\label{eq:weight_ratio}
        \frac{w(m+1)}{w(m)} = 1 + O\left(\frac{1}{\gamma_j+m}\right),
    \end{equation}
    which tends to $1$, this factor is absorbed by the decreasing factor $e^{-7\eps_1 A_j}$, yielding a ratio bounded by $e^{-6\eps_1 A_j}$ between consecutive summands in \eqref{eq:E3_interval}.

    At $m = \lceil\eps_1^{-1}\rceil$, $\phi(m) = A_j - 4\eps_1^2 A_j m^2 + O(m/\gamma_j) \le A_j - 4A_j + o(1) = -3A_j + o(1)$. Therefore \eqref{eq:E3_interval} is bounded by a geometric series with first term $O(|\Gamma(\rho_j)|e^{-3A_j})$ and ratio $O(e^{-7\eps_1 A_j})$
    \begin{equation*}
        |E_3(\mu)| \ll |\Gamma(\rho_j)|e^{-3A_j}\cdot\frac{1}{1-e^{-7\eps_1 A_j}} \ll |\Gamma(\rho_j)|e^{-3A_j}.
    \end{equation*}
    We handle all zeros with $\gamma<0$ (whether in $E_1$ or $E_3$) at once.
    For such $\rho$, we have $b_\rho^2=(\gamma_j+|\gamma|)^2\ge\gamma_j^2$ and $|a_\rho|\le1$, so
    \begin{equation*}
        \Re(\lambda_\rho)\mu \le \mu + 5\eps_1^2\alpha_j\mu - 5\eps_1^2\alpha_j\mu\gamma_j^2
        \le L_j + O(A_j) - 4\eps_1^2 A_j\gamma_j^2.
    \end{equation*}
    By Lemma~\ref{lem:slow_var}, $\alpha_j\gamma_j^2=\gamma_j^{2-o(1)}$ tends to infinity,
    hence $A_j\gamma_j^2=(\alpha_j\gamma_j^2)L_j\gg L_j$,
    and therefore $\Re(\lambda_\rho)\mu\le-2\eps_1^2A_j\gamma_j^2$ for $j$ large.
    Combining this with
    \eqref{eq:global_coeff_neg} and $N(T+1)-N(T)\ll\log T$,
    \begin{equation*}
        \sum_{\gamma<0}|\tilde c_\rho e^{\lambda_\rho\mu}|
        \ll |\Gamma(\rho_j)|\,\gamma_j^{1/2}e^{-2\eps_1^2A_j\gamma_j^2}
            \sum_{m\ge2}e^{-\pi m/2}\log m
        = o\bigl(|\Gamma(\rho_j)|e^{-2\eps_1 A_j}\bigr),
    \end{equation*}
    since $\eps_1^2A_j\gamma_j^2\gg A_j$ and $\gamma_j^{1/2}=e^{o(A_j)}$.

    Turning to $E_1$, we first consider the nearby zeros with $0 < \gamma \le 2\gamma_j$. Each term has $|e^{\lambda_\rho\mu}| \le e^{-A_j/3}$. By \eqref{eq:global_coeff_below} (for $\gamma < \gamma_j$) and \eqref{eq:global_coeff_above} (for $\gamma \ge \gamma_j$), in either case $|\tilde c_\rho| \ll |\Gamma(\rho_j)|\gamma_j^{1/2}$: when $\gamma < \gamma_j$, $(\gamma_j/\gamma)^{1/2} \le (\gamma_j/10)^{1/2} \ll \gamma_j^{1/2}$ (since $\gamma \ge 10$); when $\gamma \ge \gamma_j$, $(\gamma/\gamma_j)^{1/2} \le 2^{1/2} = O(1) \le \gamma_j^{1/2}$ since $\gamma \le 2\gamma_j$, and the exponential factor $\exp(C'|b_\rho|/\gamma_j)$ is bounded since $|b_\rho| \le \gamma_j$ gives $C'|b_\rho|/\gamma_j \le C'$. The number of such zeros is $O(\gamma_j\log\gamma_j)$ by the Riemann--von Mangoldt formula, so
    \begin{equation*}
        \sum_{\substack{\rho\in E_1\\ 0<\gamma\le 2\gamma_j}}|\tilde c_\rho e^{\lambda_\rho\mu}| \ll |\Gamma(\rho_j)|\cdot\gamma_j^{1/2}\cdot\gamma_j\log\gamma_j\cdot e^{-A_j/3} = |\Gamma(\rho_j)|\gamma_j^{3/2}\log\gamma_j\cdot e^{-A_j/3}.
    \end{equation*}
    Since $A_j/\gamma_j \to \infty$, for any sufficiently small $c > 0$ we have $\gamma_j^{3/2}\log\gamma_j\cdot e^{-A_j/3} = o(e^{-cA_j})$, in particular this is $o(|\Gamma(\rho_j)|e^{-2\eps_1 A_j})$.

    For further zeros in $E_1$ with $\gamma > 2\gamma_j$, we have $b_\rho > \gamma_j$, $a_\rho \le 0$, $a_\rho^2 \le 1$, so $\Re(\lambda_\rho) \le 5\eps_1^2\alpha_j - 5\eps_1^2\alpha_j b_\rho^2$. Since $\gamma > \gamma_j$, \eqref{eq:global_coeff_above} gives $|\tilde c_\rho| \ll |\Gamma(\rho_j)|(\gamma/\gamma_j)^{1/2}\exp(C'|b_\rho|/\gamma_j)$. Splitting the zeros into intervals of length $1$ $m < |b_\rho| \le m+1$ with $m \ge \gamma_j$ and using the weight $w(m)$ from \eqref{eq:weight_def}, we see that
    \begin{equation}\label{eq:E1_far_interval}
        \sum_{\substack{\rho\in E_1\\ \gamma > 2\gamma_j}}|\tilde c_\rho e^{\lambda_\rho\mu}| \ll |\Gamma(\rho_j)|\sum_{m\ge\gamma_j}w(m)\exp\left(\frac{C'm}{\gamma_j} + 5\eps_1^2\alpha_j\mu - 4\eps_1^2 A_j m^2\right).
    \end{equation}
    Let $\phi_1(m) := C'm/\gamma_j + 5\eps_1^2\alpha_j\mu - 4\eps_1^2 A_j m^2$ denote the exponent at index $m$. Using $\alpha_j\mu \le A_j$, we have $5\eps_1^2\alpha_j\mu \le 5\eps_1^2 A_j = O(A_j)$, so $\phi_1(m) \le C'm/\gamma_j + O(A_j) - 4\eps_1^2 A_j m^2$. The difference is $\phi_1(m+1) - \phi_1(m) = C'/\gamma_j - 4\eps_1^2 A_j(2m+1)$. For $m \ge \gamma_j$, this is $\le -4\eps_1^2 A_j(2\gamma_j + 1) + O(1/\gamma_j) \to -\infty$ as $j \to \infty$ since $A_j\gamma_j \to \infty$. At the starting point $m = \gamma_j$,
    \begin{equation*}
        \phi_1(\gamma_j) \le O(A_j) - 4\eps_1^2 A_j\gamma_j^2 = A_j\bigl(O(1) - 4\eps_1^2\gamma_j^2\bigr).
    \end{equation*}
    Since $\gamma_j \to \infty$, for $j$ large the bracketed term is $\le -2\eps_1^2\gamma_j^2$, hence $\phi_1(\gamma_j) \le -2\eps_1^2 A_j\gamma_j^2 \le -2\eps_1^2 A_j\gamma_j$ by using $\gamma_j \ge 1$. Finally, $e^{-2\eps_1^2 A_j\gamma_j}$ dominates the weight $w(\gamma_j) = O(\log\gamma_j)$, yielding that the full contribution is $o(e^{-2\eps_1 A_j})$.

    \medskip Combining all parts,
    \begin{equation}\label{eq:U_vs_E2}
        -U_j(\mu) = E_2(\mu) + O(|\Gamma(\rho_j)|e^{-2\eps_1 A_j}).
    \end{equation}

    We now apply the Tur\'an power-sum method. Recall Theorem~\ref{thm:turan} and its continuous form Corollary~\ref{cor:turan_cont} from Section~\ref{subsec:turan}.

    We first order the zeros. Enumerate the $n$ zeros of $E_2$ as $\rho^{(1)},\rho^{(2)},\ldots,\rho^{(n)}$ with $\Re(\lambda_{\rho^{(1)}}) \ge \cdots \ge \Re(\lambda_{\rho^{(n)}})$. Since $\rho_j \in E_2$ with $\lambda_{\rho_j} = 0$, we have $\Re(\lambda_{\rho^{(1)}}) \ge 0$.

    We next bound the partial sums of the coefficients required by Corollary~\ref{cor:turan_cont}. By \eqref{eq:local_coeff}, $\tilde c_{\rho^{(i)}} = \Gamma(\rho_j)(1+\theta_i)$ with $|\theta_i| \le C\eps_1^{-2}\gamma_j^{-1}$. Taking the $\Gamma(\rho_j)$ factor out, we have
    \begin{equation*}
        \sum_{i=1}^\ell\tilde c_{\rho^{(i)}} = \Gamma(\rho_j)\left(\ell + \sum_{i=1}^\ell\theta_i\right).
    \end{equation*}
    Here $n$ counts zeros with $|b_\rho| = |\gamma - \gamma_j| \le \eps_1^{-1}$, i.e., in an interval of length $2\eps_1^{-1}$ around $\gamma_j$. By the Riemann--von Mangoldt formula $N(T+1) - N(T) \ll \log T$ for the number of zeros in $[T,T+1]$, we have $n \ll \eps_1^{-1}\log\gamma_j$, so
    \begin{equation*}
        \sum_{i=1}^n|\theta_i| \le n\cdot C\eps_1^{-2}\gamma_j^{-1} \ll \eps_1^{-3}\log\gamma_j/\gamma_j \to 0.
    \end{equation*}
    For $j$ large enough that $\sum_{i=1}^n|\theta_i| < 1/2$ we then have, for every $\ell\in\{1,\ldots,n\}$,
    \begin{equation*}
        \left|\ell + \sum_{i=1}^\ell\theta_i\right| \ge \ell - \sum_{i=1}^\ell|\theta_i| \ge 1 - \tfrac{1}{2} = \tfrac{1}{2}.
    \end{equation*}
    This bound depends only on the uniform estimate $|\theta_i| \le C\eps_1^{-2}\gamma_j^{-1}$ and $n$, both independent of the ordering. Therefore
    \begin{equation}\label{eq:partial_sum}
        \min_{1\le\ell\le n}\left|\sum_{i=1}^\ell\tilde c_{\rho^{(i)}}\right| \ge \tfrac{1}{2}|\Gamma(\rho_j)|.
    \end{equation}

    We now estimate the ``loss'' factor. With $a = L_j - 6\eps_1 A_j$ and $d = \eps_1 A_j$,
    \begin{equation*}
        \frac{d}{16e(a+d)} \ge \frac{\eps_1\alpha_j}{16e},
    \end{equation*}
    where we used $a + d \le L_j = A_j/\alpha_j$. By Lemma~\ref{lem:slow_var},
    $|\log\alpha_j|=|\log\eta(\gamma_j)|=o(\log\gamma_j)$.
    Therefore the ``loss'' exponent satisfies
    \begin{equation*}
        n\log\frac{16e}{\eps_1\alpha_j} \ll \eps_1^{-1}\log\gamma_j\cdot\bigl(|\log\eps_1| + o(\log\gamma_j)\bigr) \ll \eps_1^{-1}(\log\gamma_j)^2.
    \end{equation*}
    By (H3), $\gamma_j/A_j \to 0$, so $(\log\gamma_j)^2 = o(\gamma_j) = o(A_j)$. Hence
    \begin{equation*}
        \left(\frac{d}{16e(a+d)}\right)^n \ge e^{-n\log(16e/(\eps_1\alpha_j))} \ge e^{-\eps_1 A_j}
    \end{equation*}
    for $j$ large.

    We now conclude our argument. Apply Corollary~\ref{cor:turan_cont} to $W(\mu) := E_2(\mu) = \sum_{i=1}^n\tilde c_{\rho^{(i)}}e^{\lambda_{\rho^{(i)}}\mu}$ over $\mathcal{I}_j = [a,a+d]$. Since $\Re(\lambda_{\rho^{(1)}}) \ge 0$ and $a \ge 0$, the factor $e^{a\Re(\lambda_{\rho^{(1)}})} \ge 1$. Combining \eqref{eq:partial_sum} and the ``loss'' bound gives $\mu^*\in\mathcal{I}_j$ with
    \begin{equation*}
        |E_2(\mu^*)| \ge 2\cdot 1\cdot e^{-\eps_1 A_j}\cdot\tfrac{1}{2}|\Gamma(\rho_j)| = |\Gamma(\rho_j)|e^{-\eps_1 A_j}.
    \end{equation*}
    By \eqref{eq:U_vs_E2}, the error $O(|\Gamma(\rho_j)|e^{-2\eps_1 A_j})$ is at most $\tfrac{1}{2}|\Gamma(\rho_j)|e^{-\eps_1 A_j}$ for $j$ large (since $e^{-\eps_1 A_j} \to 0$), so $|U_j(\mu^*)| \ge \tfrac{1}{2}|\Gamma(\rho_j)|e^{-\eps_1 A_j}$.

    We now derive the contradiction. Recall $D(x) := x^{-1}\int_0^x|\Delta(u)|\,du$. From the bound just obtained, $|U_j(\mu^*)| \ge \tfrac{1}{2}|\Gamma(\rho_j)|e^{-\eps_1 A_j}$, together with \eqref{eq:U_upper} evaluated at $\mu^*$,
    \begin{equation*}
        \tfrac{1}{2}|\Gamma(\rho_j)|e^{-\eps_1 A_j} \le \frac{x_j^{-\beta_j}}{\sqrt{k}}e^{22\eps_1 A_j}D(x_j) + |\Gamma(\rho_j)|e^{-c_0 A_j}.
    \end{equation*}
    By our choice of $\eps_1$, which is $\eps_1 \le \delta/48$, and $c_0 = \delta/3$, we have $c_0 - \eps_1 \ge \delta/3 - \delta/48 > \delta/4$. Hence
    \begin{equation*}
        \begin{split}
            |\Gamma(\rho_j)|e^{-c_0 A_j} &= |\Gamma(\rho_j)|e^{-\eps_1 A_j}\cdot e^{-(c_0 - \eps_1)A_j} \\
            &\le |\Gamma(\rho_j)|e^{-\eps_1 A_j}\cdot e^{-(\delta/4)A_j} \\
            &= o\bigl(|\Gamma(\rho_j)|e^{-\eps_1 A_j}\bigr)
        \end{split}
    \end{equation*}
    as $A_j$ tends to infinity, so this term can be absorbed into the left-hand side. Using $\sqrt{k} \ge 1$,
    \begin{equation}\label{eq:D_lower}
        D(x_j) \gg x_j^{\beta_j}|\Gamma(\rho_j)|e^{-23\eps_1 A_j},
    \end{equation}
    where $-\eps_1 A_j - 22\eps_1 A_j = -23\eps_1 A_j$. By the Stirling formula, $|\Gamma(\rho_j)| = e^{-B_j+O(\log\gamma_j)}$. Also $x_j^{\beta_j} = e^{\beta_j L_j} = e^{L_j}\cdot e^{-(1-\beta_j)L_j} = x_je^{-(1-\beta_j)L_j}$, and \eqref{eq:violate} gives $(1-\beta_j)L_j < (1-\delta)\alpha_j L_j = (1-\delta)A_j$. Therefore
    \begin{equation}\label{eq:D_lower2}
        D(x_j) \gg x_j\exp\bigl(-(1-\delta)A_j - B_j - 23\eps_1 A_j + O(\log\gamma_j)\bigr).
    \end{equation}

    We now derive an upper bound on $D(x_j)$ from the hypothesis. The key observation is that the function $u\mapsto ue^{-(1+\eps)\omega_\eta(u)}$ is increasing for $u \ge 1$. To see this, let $t^*(u)$ denote any minimizer in the definition of $\omega_\eta(u)$. For $1 \le u \le x$, by choosing $t = t^*(u)$ in the infimum defining $\omega_\eta(x)$,
    \begin{equation}\label{eq:omega_lipschitz}
        \omega_\eta(x) - \omega_\eta(u) \le \eta(t^*(u))\log(x/u) \le \tfrac{1}{2}\log(x/u),
    \end{equation}
    since $\eta(t) \le 1/2$. Therefore
    \begin{align*}
        ue^{-(1+\eps)\omega_\eta(u)} &= ue^{-(1+\eps)\omega_\eta(x)}\cdot e^{(1+\eps)(\omega_\eta(x)-\omega_\eta(u))} \\
        &\le ue^{-(1+\eps)\omega_\eta(x)}\cdot(x/u)^{(1+\eps)/2} \\
        &= u^{(1-\eps)/2}x^{(1+\eps)/2}e^{-(1+\eps)\omega_\eta(x)} \\
        &\le xe^{-(1+\eps)\omega_\eta(x)},
    \end{align*}
    where the last step uses $u \le x$ and $(1-\eps)/2 + (1+\eps)/2 = 1$.

    From the hypothesis $|\Delta(u)| \ll ue^{-(1+\eps)\omega_\eta(u)}$ for $u \ge u_0$ and the above monotonicity, $|\Delta(u)| \ll xe^{-(1+\eps)\omega_\eta(x)}$ for all $u \in [u_0, x]$. Splitting the integral at $u_0$,
    \begin{equation*}
        D(x) = \frac{1}{x}\int_0^{u_0}|\Delta(u)|\,du + \frac{1}{x}\int_{u_0}^x|\Delta(u)|\,du \le \frac{O(u_0)}{x} + \frac{x-u_0}{x}\cdot xe^{-(1+\eps)\omega_\eta(x)}.
    \end{equation*}
    The second term is $\le xe^{-(1+\eps)\omega_\eta(x)}$. By \eqref{eq:omega_lipschitz} applied with $u = 2$, $\omega_\eta(x) \le \omega_\eta(2) + \tfrac{1}{2}\log(x/2) \le \tfrac{1}{2}\log x + O(1)$, so $xe^{-(1+\eps)\omega_\eta(x)} \ge cx^{(1-\eps)/2}$. This tends to infinity, so for $x$ large enough it exceeds $1$, and therefore the first term $O(u_0/x) = o(1)$ is bounded by $xe^{-(1+\eps)\omega_\eta(x)}$ as well. Combining all these together we have
    \begin{equation}\label{eq:D_from_hyp}
        D(x_j) \ll x_je^{-(1+\eps)\omega_\eta(x_j)}.
    \end{equation}
    By Lemma~\ref{lem:localization},
    \begin{equation}\label{eq:D_upper}
        D(x_j) \ll x_j\exp\bigl(-(1+\eps)(A_j + B_j) + o(A_j)\bigr).
    \end{equation}

    Comparing \eqref{eq:D_lower2} and \eqref{eq:D_upper}, the lower bound forces
    \begin{equation*}
        -(1-\delta)A_j - B_j - 23\eps_1 A_j + O(\log\gamma_j) \le -(1+\eps)(A_j + B_j) + o(A_j),
    \end{equation*}
    which means
    \begin{equation*}
        (1+\eps)A_j + (1+\eps)B_j - (1-\delta)A_j - B_j - 23\eps_1 A_j \le o(A_j) + O(\log\gamma_j),
    \end{equation*}
    i.e., $(\delta + \eps - 23\eps_1)A_j + \eps B_j \le o(A_j)$ by absorbing $O(\log\gamma_j) = o(A_j)$. Recall that we chose $\eps_1 \le \eps/48$, so $23\eps_1 \le 23\eps/48 < \eps/2$, and hence $\delta + \eps - 23\eps_1 > \delta + \eps/2$. Since $B_j \ge 0$, we conclude $(\delta + \eps/2)A_j \le o(A_j)$, contradicting $\delta > 0$ and $A_j \to \infty$.
\end{proof}

\subsection{Proof of Theorem~\ref{thm:main2}}

\begin{proof}[Proof of Theorem~\ref{thm:main2}]
    We apply Proposition~\ref{prop:Gk} with $\varphi = (1\pm\eps)\omega_\eta$. We verify both hypotheses. Firstly, $\omega_\eta(x) \to \infty$ as $x \to \infty$ since $\eta(t) > 0$ on $[2,\infty)$, hence $\varphi(x) \to \infty$. Secondly, since $\eta(t) \le 1/2$ gives $\omega_\eta(x) \le \frac{1}{2}\log x + O(1)$, we have $(1+\eps)\omega_\eta(x) \le \frac{1+\eps}{2}\log x + O(1)$, so
    \begin{equation*}
        xe^{-(1+\eps)\omega_\eta(x)} \ge c\,x^{(1-\eps)/2} \to \infty.
    \end{equation*}
    The case $\varphi = (1-\eps)\omega_\eta$ can be proved similarly, even easier. Proposition~\ref{prop:Gk} therefore reduces both directions to Theorem~\ref{thm:main1}.
\end{proof}

\section{Proof of the Upper Bound for \texorpdfstring{$W(x)$}{W(x)}}
\label{sec:upper_bound}

We prove that $W(x) \ll_\eps xe^{-(1-\eps)\omega(x)}$ for any $\eps > 0$.

\begin{proof}
Using $|\Gamma(\rho+1)| = |\rho|\cdot|\Gamma(\rho)| \le \sqrt{\beta^2+\gamma^2}\cdot|\Gamma(\rho)| \le \sqrt 2\gamma|\Gamma(\rho)|$ for $0 < \beta < 1$, $\gamma > 1$,
\begin{equation*}
    W(x) = \sum_{0<\gamma\le x}\frac{|\Gamma(\rho+1)|x^\beta}{\gamma} \le \sqrt 2\sum_{0<\gamma\le x}|\Gamma(\rho)|x^\beta \le 2\sum_{\gamma>0}|\Gamma(\rho)|x^\beta,
\end{equation*}
where zeros with $\gamma > x$ contribute $O(e^{-x})$ by the exponential decreasing of $|\Gamma(\rho)|$.

By definition of $\omega(x)$, every zero satisfies $|\Gamma(\rho)|x^\beta \le xe^{-\omega(x)}$. For any $\eps\in(0,1)$, split each summand as
\begin{equation*}
    |\Gamma(\rho)|x^\beta = \bigl(|\Gamma(\rho)|x^\beta\bigr)^{1-\eps}\cdot|\Gamma(\rho)|^\eps\cdot x^{\eps\beta}.
\end{equation*}
The first factor is $\le(xe^{-\omega(x)})^{1-\eps} = x^{1-\eps}e^{-(1-\eps)\omega(x)}$. Combining with $x^{\eps\beta}$ gives $x^{1-\eps+\eps\beta} = x^{1-\eps(1-\beta)} \le x$ since $\beta\le 1$. Therefore
\begin{equation*}
    |\Gamma(\rho)|x^\beta \le xe^{-(1-\eps)\omega(x)}\cdot|\Gamma(\rho)|^\eps.
\end{equation*}
By the Stirling formula, $|\Gamma(\rho)|^\eps \ll \gamma^{\eps/2}e^{-\eps\pi\gamma/2}$. Applying the Riemann--von Mangoldt formula we have
\begin{equation*}
    \sum_{\gamma>0}|\Gamma(\rho)|^\eps \ll_\eps \int_{10}^\infty t^{\eps/2}e^{-\eps\pi t/2}\log t\,dt < \infty.
\end{equation*}
Therefore $W(x) \ll_\eps xe^{-(1-\eps)\omega(x)}$. 
\end{proof}

\section*{Acknowledgments}
The author would like to thank his supervisor, Professor A. I. Suriajaya, for her kind guidance and support.
The author used Codex to assist with drafting and editing the manuscript, checking mathematical arguments, and revising parts of the proofs.
The tool also suggested organizing parts of the proofs into separate lemmas to improve readability.
The author takes full responsibility for the content and mathematical correctness of this paper.

\nocite{*}
\bibliographystyle{amsplain}
\bibliography{b}

\end{document}